\documentclass[12pt]{amsart}
\long\def\eatit#1{}

\usepackage{pgf,tikz,hyperref,comment}
\usepackage{amsmath, graphicx}
\usepackage{amssymb,color}
\usepackage{enumitem}

\newtheorem{theorem}{Theorem}[section]
\newtheorem{proposition}[theorem]{Proposition}
\newtheorem{lemma}[theorem]{Lemma}
\newtheorem{corollary}[theorem]{Corollary}

\newtheorem*{theorem*}{Theorem}

\newtheorem{assump}{Theorem}    
\newenvironment{myassump}[2][]
  {\renewcommand\theassump{#2}\begin{assump}[#1]}
  {\end{assump}}

\theoremstyle{definition}
\newtheorem{problem}[theorem]{Problem}
\newtheorem{prob*}{Problem}
\newtheorem{question}[theorem]{Question}
\newtheorem{definition}[theorem]{Definition}

\newtheorem{example}[theorem]{Example}
\newtheorem{remark}[theorem]{Remark}

\newtheorem{thevarthm}[theorem]{\varthmname}

\newenvironment{varthm*}[1]{\trivlist\item[]{\bf #1.}\it}{\endtrivlist}

\newcommand{\PP}{ \ensuremath{\mathbb{P}}}

\begin{document}

\author[L.~Chiantini]{Luca Chiantini}
\address[L.~Chiantini]{Dipartmento di Ingegneria dell'Informazione e Scienze Matematiche, Universit\`a di Siena, Italy}
\email{luca.chiantini@unisi.it}

\author[{\L}.~Farnik]{{\L}ucja Farnik}
\address[{\L}.~Farnik]{Department of Mathematics, University of the National Education Commission, Krakow,
   Podcho\-r\c a\.zych~2,
   PL-30-084 Krak\'ow, Poland}
\email{lucja.farnik@gmail.com}

\author[G.~Favacchio]{Giuseppe Favacchio}
\address[G.~Favacchio]{Dipartimento di Ingegneria, Universit\`a degli studi di Palermo,
Viale delle Scienze,  90128 Palermo, Italy}
\email{giuseppe.favacchio@unipa.it}

\author[B.~Harbourne]{Brian Harbourne}
\address[B.~Harbourne]{Department of Mathematics,
University of Nebraska,
Lincoln, NE 68588-0130 USA}
\email{brianharbourne@unl.edu}

\author[J.~Migliore]{Juan Migliore} 
\address[J.~Migliore]{Department of Mathematics,
University of Notre Dame,
Notre Dame, IN 46556 USA}
\email{migliore.1@nd.edu}

\author[T.~Szemberg]{Tomasz Szemberg}
\address[T.~Szemberg]{Department of Mathematics, University of the National Education Commission, Krakow,
   Podcho\-r\c a\.zych~2,
   PL-30-084 Krak\'ow, Poland}
\email{tomasz.szemberg@gmail.com}

\author[J.~Szpond]{Justyna Szpond}
\address[J.~Szpond]{Department of Mathematics, University of the National Education Commission, Krakow,
   Podcho\-r\c a\.zych~2,
   PL-30-084 Krak\'ow, Poland}
\email{szpond@gmail.com}

\thanks{Chiantini and Favacchio are members of INdAM - GNSAGA}
\thanks{Farnik was partially supported by the National Science Centre, Poland, grant 2018/28/C/ST1/00339.
}
\thanks{Favacchio was partially supported by “Piano straordinario per il miglioramento della qualità della ricerca e dei risultati della VQR
2020-2024 - Misura A” of the Università degli studi di Palermo}
\thanks{Harbourne was partially supported by Simons Foundation grant \#524858.}
\thanks{Migliore was partially supported by Simons Foundation grant \#839618.
}
\thanks{Szemberg and Szpond were partially supported by the National Science Centre, Poland, grant 2019/35/B/ST1/00723.}


\begin{abstract}
In this note we introduce the notion of  $(b,d)$-geprofi sets and study their basic properties. These are sets of $bd$ points in $\PP^4$ whose projection from a general point to a hyperplane is a full intersection, i.e., the intersection of a curve of degree $b$ and a surface of degree $d$. 
We show that such nontrivial sets exist if and only if $b\geq 4$ and $d\geq 2$.

Somewhat surprisingly, for infinitely many values of $b$ and $d$ there exist such sets in  linear general position. The note contains open questions and problems.
\end{abstract}


\title{Finite sets of points in $\mathbb P^4$ with special projection properties   }

\date{July 1, 2024}
\maketitle

\markleft{L.~Chiantini, \L.~Farnik,  G.~Favacchio, B.~Harbourne, J.~Migliore, T.~Szemberg and J.~Szpond}

\section{Introduction}

We work over the field $\mathbb C$ of complex numbers. For a finite set of points $Z \subset \PP^N$ (i.e., a $0$-dimensional reduced scheme) we write $H_Z$ for the Hilbert function and $h_Z$ for the corresponding $h$-vector.

In recent years considerable attention has been given to \emph{geproci sets} of points in projective spaces, 
i.e., finite sets of points, whose projection from a general point to a hyperplane is a complete intersection. Families of examples of such sets in $\PP^3$ have been constructed and extensively studied in \cite{POLITUS1}, \cite{POLITUS2}, \cite{POLITUS3}, \cite{cortona} and \cite{POLITUS4b}. Currently no examples of geproci sets in higher dimensional spaces are known and there is  accumulating evidence to the effect that possibly no such sets  exist at all.

In this paper we obtain an interesting theory in higher dimensions by modifying the definition of geproci sets in such a way
that for $\PP^3$ we get back the original definition of geproci sets.
Note that on $\PP^2$ a curve is also a hypersurface, hence a complete intersection is the intersection of a curve and a hypersurface. So instead of asking for the projection to a hyperplane to be a complete intersection we introduce a new concept, that of a \emph{geprofi} set. On $\PP^4$, we do this by asking that the projection be the intersection of a curve and a surface in the hyperplane.
Making this change requires the development of entirely new techniques.

To be even more general, we begin with the following definition.

\begin{definition}[Full intersection]\label{def: full inter}
   We say that a reduced finite set of points $Y$ in a projective space $\PP^M$ is a {\it full intersection} if there exist varieties $U,V\subset \PP^M$ of complementary dimensions (i.e., $\dim(U)+\dim(V)=M$), such that
   \begin{itemize}
       \item $U\cap V=Y$ and
       \item $\# Y =\deg(U)\cdot\deg(V)$.
   \end{itemize}
\end{definition} 

Of course, if $U$ and $V$ are complete intersections, then so is $Y$. In particular, if $M=2$ and $U$ and $V$ are planar curves, then $Y$ is a complete intersection. On the other hand there are many examples of full intersections which are not complete intersections, so it is clearly a much wider class of subschemes. For example, a reduced hypersurface section of a curve is a full intersection, but it is a complete intersection if and only if the curve is a complete intersection.

\begin{definition}[Geprofi property]\label{def: geprofi}  
We say that a set of points $Z\subset\PP^N$ has the \emph{geprofi property} (general projection full intersection) if there exists $M\geq 2$ such that for the projection $\pi_{\Lambda}:\PP^N\dashrightarrow\PP^M$
from a \textbf{general} linear subspace $\Lambda$ of $\PP^N$ of dimension $(N-M-1)$, the image $Y=\pi_{\Lambda}(Z)$ is a full intersection.
\end{definition}      

In this work we restrict our attention to geprofi sets of points in $\PP^4$ and their projections to $\PP^3$, even though such sets exist in projective spaces of higher dimension. It is convenient to work here with the following definition.

\begin{definition}[A $(b,d)$-geprofi set]
We say that a  finite set of points $Z\subset\PP^4$ is a $(b,d)$-\emph{geprofi set} if, for the projection $\pi_{P}:\PP^4\dashrightarrow\PP^3$
from a \textbf{general} point $P\in \PP^4$, the image $\pi_{P}(Z)$ is a full intersection of a curve of degree $b$ and a surface of degree $d$.
\end{definition}      

 Suppose $Z \subset \PP^4$ is a $(b,d)$-geprofi set. There are two
 fundamental kinds of such sets.
\begin{itemize}[leftmargin=45pt,]
    \item[Type 1:] $Z$  is already  the intersection of a curve of degree $b$ and a surface of degree $d$ in $\PP^4$. (We will also refer to sets $Z$ of type 1 as {\it trivial $(b,d)$-geprofi sets}, since it is obvious that the general projection will have the same property.)
    \item[Type 2:] $Z$  is not the intersection of a surface of degree $d$ and a curve of degree $b$ in $\PP^4$. (We will refer to such sets as {\it nontrivial $(b,d)$-geprofi sets}.)
\end{itemize}
Contrary to geproci sets, it may happen that the numbers $b$ and $d$ are not uniquely assigned to a set of points $Z$. In fact, it may even happen that the same set of points $Z$ is of type 1 for the pair of integers $(b,d)$, while it is of type 2 for some other pair of integers $(b',d')$; see Example \ref{ex: 3 4 and 6 2 example}.

We give a careful description of  geprofi sets of type 1 in Section \ref{sec: triv} -- see Theorem \ref{classi} and Proposition \ref{one point} for a description of their structure.
The rest of the paper is devoted to geprofi sets of type 2; 
indeed, the study of properties of $Z$ that force the general projection to be a full intersection in $\PP^3$ even though  $Z$ itself is not the intersection in $\PP^4$ of a curve of degree $b$ and a surface of degree $d$ seems more interesting.

This situation is analogous to $(a,b)$-geproci sets in $\PP^3$: either $Z$ is the intersection of curves of degrees $a$ and $b$ in $\PP^3$, or it is not, and the general projection is a complete intersection for subtle reasons. In \cite[Proposition 3.1]{CM}, the work of Diaz \cite{diaz} and Giuffrida \cite{giuffrida} was extended to show that when $|Z| = ab$ and $Z$ is the intersection of a curve $C_a$ of degree $a \geq 2$ and a curve $C_b$ of degree $b \geq a$ in $\PP^3$, then both $C_a$ and $C_b$ are necessarily unions of disjoint lines. If $a>2$, $C_a$ and $C_b$ lie on a common smooth quadric surface in opposite rulings, forming what are called {\it grids}. In the geprofi context, we will show that the property that $Z$ is the intersection of a curve of degree $b$ and a surface of degree $d$ in $\PP^4$, with $|Z| = bd$, is restrictive but not as much as in the geproci situation. The heart of this paper studies geprofi sets of type 2.

In \cite[Question 0.2]{POLITUS1} the authors asked about the existence of {\it geproci} sets in $\PP^3$ in linear general position (henceforth denoted LGP), apart from 4 general points. This question is still open. On the other hand, there are infinitely many geprofi sets in $\PP^4$  in LGP, and in this paper we also begin the study of sets with this property.

There are three classification problems that are of interest to us. 

\begin{enumerate}

\item For which $b,d$ does there exist a $(b,d)$-geprofi set in $\PP^4$? (We have a complete answer in Theorem \ref{numerical classification}.)

\item \label{2nd} For which $b,d$ does there exist a $(b,d)$-geprofi set in LGP in $\PP^4$? (We have a partial answer in Theorem \ref{2b LGP class} and Remark \ref{LGP table}.) 

\item \label{3rd} For which $b,d$ does there exist a $(b,d)$-geprofi set in $\PP^4$ such that the image under a general projection is the full intersection of an irreducible curve of degree $b$ and an irreducible surface of degree $d$? (We have partial answers, e.g., Corollary \ref{find geprofi by curves}.)
\end{enumerate}

More precisely, the following summarize the main results of this paper in the direction of classifying geprofi sets. 

\begin{myassump}{A} {\em [Theorem \ref{numerical classification}]} \label{thm: A}
A nondegenerate, nontrivial $(b,d)$-geprofi set exists in $\PP^4$ if and only if  $b \geq 4$ and $d \geq 2$.  
\end{myassump}

The next two results concern sets of points in LGP.

\begin{myassump}{B} {\em [Theorem \ref{2b LGP class}] } \label{thm: B}
  Let $Z$ be a set of $2b$ points in LGP in $\PP^4$. 
 Then  $Z$ is a $(b,2)$-geprofi set if and only if one of the following occurs.
    \begin{itemize}
\item[(i)]$b=3$ (in fact such $Z$ is always a trivial geprofi set).

\item[(ii)] $b=4$. 

\item[(iii)] $b=5$, $Z$ is arithmetically Gorenstein, and either it lies on a rational quartic curve, or it is linked by a system of quadrics to a degenerate scheme of length $6$.

\item[(iv)]$b\geq 6$ and $Z$ lies on a rational normal curve.  (If $b$ is odd then $Z$ is arithmetically Gorenstein.)
    \end{itemize}
\end{myassump}

\begin{myassump}{C} {\em [cf. Corollary \ref{find geprofi by curves}]} \label{thm: C}
 For fixed $b$, assume $d > -\frac{10}{b}+2$ and $b(d-1) + 5 \leq \binom{d+3}{3}$. Then a nontrivial $(b,d)$-geprofi set $Z$ in LGP exists on a smooth curve of degree $b$ in $\PP^4$. In particular, for fixed $b$ and all $d \gg 0$, a $(b,d)$-geprofi set in LGP exists.
\end{myassump}

To prove this last result we combine a theorem of Gruson-Lazarsfeld-Peskine \cite{GLP} with  results of Ballico-Ellia \cite{BE-inv}. In the projection, the curve comes from $\PP^4$, but the surface appears only after projecting.

We remark that we do not know what happens in the other direction. Indeed, our work suggests the following two interesting open questions.

\begin{question}

\begin{enumerate}

\item For fixed $d$ and $b \gg 0$, does there exist a $(b,d)$-geprofi set in LGP? (See Question \ref{final question} for more in this direction.)
    
 \item  Does there exist a set of points $Z \subset \PP^4$ such that $Z$ is $(b,d)$-geprofi, but yet $Z$ does not lie on a curve of degree $b$ in $\PP^4$?
 (This question is not limited to LGP sets.)      
    \end{enumerate}

\end{question}

 In addition to the above, throughout the paper we give many constructions that produce geprofi sets in interesting ways, and we use some of them toward proofs of the above main results. We also give connections to Weddle loci (cf. \cite[Chapter 2]{POLITUS1}) and to Lefschetz properties and unexpected hypersurfaces \cite{HMNT}.

 Although the notion of geprofi sets was motivated by that of geproci sets, this paper can be read independently from the existing literature on geproci sets. The techniques introduced here are quite different from those used to study geproci sets.


\section{Geprofi sets of type 1} \label{sec: triv}

In the discussion  following Definition  \ref{def: geprofi}, we noted that the  geprofi property for a set $Z \subset \PP^4$ is obvious if $Z$ is already the intersection of a curve of degree $b$ and a surface of degree $d$. In this section we give some idea of how restrictive this condition is, and the place to start is to compare with the situation of geproci sets in $\PP^3$. 

\begin{definition}
A finite set of points $Z\subset\PP^3$ is an {\it $(a,b)$-geproci set} if its projection from a general point in $\PP^3$ to a plane $H\subset \PP^3$ is the transverse intersection of curves in $H$ of degrees $a$ and $b$. 
\end{definition}

The question  of the existence of geproci sets was raised by Polizzi \cite{Polizzi}. It was observed by Panov (see \cite{Polizzi}) that there is an obvious way  to find nondegenerate geproci sets:  if $Z \subset \PP^3$ consists of the intersection  points of a set  $C$ of $a$ lines in one ruling of a smooth quadric $\mathcal{Q}$ and  a set $C'$ of $b$ lines in the other ruling, then $Z$ is clearly $(a,b)$-geproci. This is because then a general projection   $\pi(Z)$ is the complete intersection  of $\pi(C)$ and $\pi(C')$.

In fact, this is essentially  the only such example. In \cite[Proposition 3.1]{CM}, the (a priori more general)  situation  was studied where $Z$ is the intersection  of a curve $C$ of degree $a$ and a curve  $C'$ of degree $b$, with $2 \leq a \leq b$. It was shown that if $Z$ is not  degenerate, and if $a \geq 3$, then in fact $Z$ can only arise in the setting described by Panov. If $a=2$, both  $C$ and $C'$ must consist of pairwise skew lines but $C'$ need not lie on a quadric. In short, $Z$ must be a grid. This was extended in \cite[Theorem A]{POLITUS3}, where it was shown that if $Z$ is $(a,b)$-geproci and lies on a curve $C$ of degree $c \leq \max\{a,b\}$, where $c$ is the least degree among curves containing $Z$, then  $C$ must be a set of skew lines. In \cite{POLITUS1}, such a geproci set, when not a grid, is called a {\it half grid}. In \cite{POLITUS2}, \cite{POLITUS3} and \cite{cortona}, the geometry  and combinatorics of half grids was studied. Their geometry  is very  striking and nontrivial, but nevertheless the assumption about the curve of degree $c$  clearly  turns out to be  a very  restrictive one. 

Turning  now to geprofi sets in $\PP^4$, one would  like to know to what extent such behavior continues to hold. We first make a basic definition.

\begin{definition}
A  set of points $Z \subset \PP^4$ is {\it a trivial $(b,d)$-geprofi set} if $|Z| = bd$ and $Z$ is the intersection  of a  curve $C$ of degree $b$ and a  surface $S$ of degree $d$. 
\end{definition}

\begin{remark}
In the above definition we do not make any assumption that $Z$ is nondegenerate. If it is,  this forces both $C$ and $S$ to be nondegenerate. But what if $Z$ is degenerate? There are (a priori) several possibilities.

\begin{enumerate}

\item Both $C$ and $S$ are nondegenerate. 

\item $C$ is degenerate and $S$ is nondegenerate.

\item $S$ is degenerate and $C$ is nondegenerate.

\begin{enumerate}
    \item $S$ is a plane.

    \item $S$ spans a $\PP^3$.
\end{enumerate}

\item Both $C$ and $S$ are degenerate.

\begin{enumerate}
    \item $C$ and $S$ live in the same $\PP^3$.

    \item $C$ and $S$ do not live in the same $\PP^3$.
\end{enumerate}
    
\end{enumerate}

\noindent The interesting cases are those where $C \cup S$ is nondegenerate in $\PP^4$, i.e., all except (4a).

\end{remark}

To begin with it is very interesting to observe that, contrary to the geproci situation, the same set of points $Z$ can be geprofi for various values of $b$ and $d$.

\begin{example}[A double geprofi set of points]\label{ex: 3 4 and 6 2 example}
    Let $Z$ be a set of $12$ points distributed evenly on $3$ general lines in $\PP^4$. Then $Z$ is a trivial $(3,4)$-geprofi set. Indeed, let $\pi$ be a general projection from $\PP^4$ to $\PP^3$. Then one can take $C$ to be the union  of the 3 lines, and $S$ to be the three planes spanned by  subsets of 3 of the 12 points (there is obviously more than  one way to do  this), so $Z = C \cap S$.
    
    On the other hand, the set $Z$ is also $(6,2)$-geprofi. To see this, consider the unique quadric $\mathcal Q$ containing the images of the lines $\pi$ as a surface of degree $2$, and the union of $6$ lines joining pairs of points in the projection of $Z$ (and not contained in $\mathcal Q$) as a curve of degree $6$. Note that $Z$ is not a {\it trivial} $(6,2)$-geprofi set. Indeed  $Z$ is not contained in quadric surfaces in $\PP^4$ because it is nondegenerate and cannot be contained in a pair of disjoint planes.
\end{example}

Now we give a simple example where the curve does not come from $\PP^4$. 

\begin{example} \label{six pts}
    Let $Z \subset \PP^4$ be a set of 6 general points. It is easy to see that $Z$ is a trivial $(3,2)$-geprofi set using 2 planes and 3 lines (in more than one way). On the other hand, a general projection of $Z$ lies on a twisted cubic curve $C$ and on a smooth quadric surface $\mathcal Q$ not containing $C$. Notice that neither the twisted cubic $C$ nor the quadric $\mathcal Q$ is the projection from $\PP^4$, while the 3 lines in $\PP^3$ and the planes are projected from $\PP^4$.
    
     If $Z \subset \PP^4$ consists of a general set of $b \geq 4$ points on each of two general planes, it is clear  (by considering  lines spanned by  pairs of points) that $Z$ is a trivial $(b,2)$-geprofi set. But we can also view the  projection  as a full intersection where the  curve does {\it not} come as a projection from $\PP^4$:   the general projection is the full intersection of two planes in $\PP^3$ with a curve of degree $b$ that is not a union of lines, by taking 6 points at a time (3 on each plane) and using a twisted cubic for those as above, and once there remain fewer than 6 points then using lines to cut out pairs of points. 
    \end{example}

We next consider a special class of trivial geprofi sets that is modeled after grids in $\PP^3$.

\begin{definition} (Hypergrid).\label{def: hypergrid}
A $(b,d)$-geprofi set $Z \subset \PP^4$ is a {\it $(b,d)$-hypergrid} if there is a set of $b$  skew  lines and a set of $d$ planes such that the union of any  two planes is nondegenerate, and such that each line meets each plane in exactly one point, and these $bd$ points of intersection are the  points of $Z$.
\end{definition}

The next result shows that hypergrids exist for arbitrary values of $b$ and $d$.

\begin{lemma}[Existence of hypergrids]
There exist $(b,d)$-hypergrids for all $b,d\geq 1$.
\end{lemma}

\begin{proof}
Let $\alpha$ be the Segre embedding  $\alpha:\PP^1\times\PP^2 \longrightarrow\PP^5$ and let $Z'$ be the image of $D\times B\subset \PP^1\times \PP^2$ under $\alpha$, where $D\subset\PP^1$ is $d$ distinct points and $B\subset\PP^2$ is a set of $b$ distinct points. Let $\pi:\PP^5\longrightarrow\PP^4$ be a general projection. Then $Z=\pi(Z')$ is a $(b,d)$-hypergrid.
\end{proof}

Trivial $(b,d)$-geprofi sets do not need to be hypergrids. The next three examples show the extent to which this is true, even if the assumption  of being a trivial $(b,d)$-geprofi set is still quite restrictive (an assertion that we will make precise in this section).

\begin{example} \label{triv ex 1}
($C$ does not need to be a union of lines.)

Let 

\begin{itemize} 
    \item $\Lambda$ be a general plane in $\PP^4$.
    \item $\lambda$ be a general line in $\PP^4$. 
    \item $L_1,\dots,L_d$ be general lines in the plane $\Lambda$.
    \item $P_1,\dots,P_d$ be general points on the line  $\lambda$.
    \item $A_i$ be the plane spanned by $L_i$ and $P_i$ for $1 \leq i \leq d$.
    \item $S = A_1 \cup \dots \cup A_d$.
    \item $B$ be a  curve on $\Lambda$ of degree $b-1$ meeting  $(L_1 \cup \dots \cup L_d)$ in $d(b-1)$ distinct points.
    \item $C = B \cup \lambda$.
\end{itemize}

\noindent Then $|S \cap C| = bd$. In this example note the following properties.

\begin{enumerate}
\item $C$ and $S$ are both nondegenerate in $\PP^4$.

\item $C$ is not necessarily a union of lines.

\item $S$ is a union of planes.

\item Every component of $C$ spans at most a plane.

\item If $S = \cup S_i$ and $C = \cup C_j$ are the irreducible decompositions then for each $i$ and $j$, $S_i \cup C_j$ spans  a hyperplane.
\end{enumerate}
\end{example}

\begin{example} \label{triv ex 2}
($S$ does not need to be a union of planes and $C$ does not need to be a union of lines, but $C$ lies on a plane.)

Let $\Lambda$ be a plane in $\PP^4$ and let $H_1,\dots,H_r$ be hyperplanes  containing $\Lambda$. For $1 \leq i \leq r$ let $S_i$ be a general surface of degree $d_i$ in $H_i$, where $d_1 + \dots + d_r = d$, and let $S = \cup S_i$. Let $C = \cup C_i$ be a union of general curves of degrees $b_1,\dots,b_s$ in $\Lambda$, where $b_1 + \dots + b_s = b$. Note that $\Lambda \cap S_i$ is a plane curve of degree $d_i$ in $\Lambda$. Thus $|C \cap S| = (b_1 + \dots + b_s) (d_1 + \dots + d_r) = bd$. In this example note the following properties.

\begin{enumerate}
\item $S$ is nondegenerate in $\PP^4$ but $C$ spans only a plane.

\item Each component of $S$ spans a hyperplane (in particular it is degenerate).

\item Neither $S$ nor $C$ is a union of linear varieties (except when $d_i = 1$ or $b_j = 1$).

\item Any union $C_j \cup S_i$ spans a hyperplane (specifically $H_i$).
\end{enumerate}
\end{example}

\begin{example} \label{triv ex 3}
(The components of $C$ can  each be a curve that spans a hyperplane.)

Let $S$ be a plane in $\PP^4$, so $\deg S = d = 1$. Fix $r, b \geq 1$ and fix integers $b_1,\dots,b_r$ with $b_1 + \dots + b_r = b$. For $1 \leq j \leq r$, let $H_j$  be a general hyperplane containing $S$ and let $C_j$ be a nondegenerate curve in $H_j$ of degree $b_j$ meeting $S$ in $b_j$ distinct points. Let $C = C_1 \cup \dots \cup C_r$ with $\sum \deg(C_i) = b$. Then $|C \cap S| = bd = b$.  

In this example note the following properties.

\begin{enumerate}

\item $S$ is a plane.

\item $C$ is nondegenerate.

\item Each $C_i$ spans a $\PP^3$.

\item For each $j$, $S \cup C_j$ spans the hyperplane $H_j$. 

\end{enumerate}
\end{example}

\begin{example} \label{2planecurves}
(Both $C$ and $S$ can be nondegenerate.)
Let $P$ be a point in $\PP^4$ and let $\lambda_1,\dots,\lambda_4$ be general lines through $P$. Let
\[
\begin{array}{lcl}
\Pi_1 & = & \hbox{ plane spanned by } \lambda_1, \lambda_2, \\
\Pi_2 & = & \hbox{ plane spanned by } \lambda_3, \lambda_4.
\end{array}
\]
Note that $\Pi_1$ and $\Pi_2$ meet only  at $P$ so they span $\PP^4$.
For $i = 1,2$ let $B_i$ be a general plane curve in $\Pi_i$ of degree $b_i$.
Let
\[
\begin{array}{lcl}
\pi_1 & = & \hbox{ plane spanned by } \lambda_1, \lambda_3, \\
\pi_2 & = & \hbox{ plane spanned by } \lambda_2, \lambda_4.
\end{array}
\]
Let $L_1,\dots, L_q$ be lines joining generally  chosen points of $\pi_1$ to generally chosen points of $\pi_2$.

Let $S = \pi_1 \cup \pi_2$ and let $C = B_1 \cup B_2 \cup L_1 \cup \dots \cup L_q$. Let $Z = C \cap S$. Then $|Z| = |S \cap C| = (\deg S)(\deg C)$ and so $Z$ is a trivial and nondegenerate geprofi set.

 In this example note the following properties.

\begin{enumerate}

\item $S$ is the union  of two planes.

\item $C$ and $S$ are both nondegenerate, as is $C \cap S$.

\item Each component of $C$ is either a line or a plane curve. There are arbitrarily  many  lines, but only  two plane curves that are not lines.

\item The union of any  component of $S$ and any  component of $C$ spans a hyperplane.

\end{enumerate}

\end{example}

As was noted above, in $\PP^3$ the (nondegenerate) ``trivial geproci sets" in $\PP^3$ turn out to be exactly the grids. In the above examples  we see that in contrast, trivial  geprofi sets in $\PP^4$ do not need to be the intersection of a union of planes with a union of lines, although such sets do exist. However, our examples are ``not far" from this and it is worth seeing what properties must be true for trivial geprofi sets.

\begin{remark} \label{reduce to irred}
If $C = C_1 \cup \dots \cup C_r$ and $S = S_1 \cup \dots \cup S_t$ form a trivial geprofi set (i.e., $| C \cap S| = (\deg C)(\deg S)$) then the same is true of any individual $C_i$ with any individual $S_j$.  So as a first step it makes sense to assume that $C$ and $S$ are irreducible.
\end{remark}

\begin{proposition} \label{triv irred}
Let $C,S \subset \PP^4$ be a reduced and irreducible curve and surface respectively with $\deg C = b$ and $\deg S = d$. Assume that  $C \not \subset S$. Then $|C \cap S| = bd$ if and only  if $C$ and $S$ are both degenerate, lying in the same hyperplane, and the scheme-theoretic intersection $C \cap S$ is reduced.
\end{proposition}

\begin{proof}
Our assumptions imply that $|C \cap S|$ is finite, since $C$ is irreducible and $C \not \subset S$. Thus if $C$ and $S$ lie inside the same hyperplane and their intersection is reduced, B\'ezout's theorem gives $|C \cap S| = bd$. 

For the converse we assume that $|C \cap S| = bd$ and we will break into cases.

\begin{itemize}

\item Assume $d = 1$ (so $S$ is a plane) and $b$ is arbitrary.
Choose a general point $P$ of $C$ and let $H$ be the hyperplane spanned by $S$ and $P$. We have assumed that $S$ contains exactly $(b)(1) = b$ points of $C$, so $H$ contains at least $b+1$ points of $C$. Then by B\'ezout's theorem and the irreducibility  of  $C$, $H$ contains $C$. By  construction $H$ also contains $S$, so we are done.

\vspace{.1in}

\item Assume that $d \geq 2$ but that $S$ is degenerate. As before assume $b$ is arbitrary. Since $S$ is not a plane, $S$ must span $H = \PP^3$. If $C \subset H$ we are done, so assume otherwise. Then $| C \cap S | \leq |C \cap H | \leq b < bd$ and we have a contradiction.

\end{itemize}

\noindent At this point we have finished the situation where $S$ is degenerate, and we assume from now on that $S$ is nondegenerate. In particular, $d \geq 3$. We want to show that it cannot happen that $|S \cap C| = bd$.

\begin{itemize}

\item 
Assume that $C$ spans at least $\PP^3$. (In particular, $b \geq 3$.) We are assuming that we have a nondegenerate, irreducible surface $S$ of degree $d$ and an irreducible curve $C$ of degree $b \geq 3$ such that $|C \cap S| = bd$. Let $H$ be a general hyperplane. First, let $P_1$ be any point of $\PP^4 \backslash S$. Since the cone over $S$ with vertex $P_1$ is a 3-fold, a general line through $P_1$ avoids $S$. In particular, if $P_1 \in C$ is a general point of $C$ (so $P_1 \notin S$) and $P$ is a general point in $\PP^4$ then the line $\overline{PP_1}$ is disjoint from $S$. Denote by $\pi_P$ the projection from $P$ to $H$.
Thus
\[
\{ P \in \PP^4 \ | \pi_P(C) \not \subset \pi_P(S) \}
\]
contains an open set, since $\pi_P(P_1) \notin \pi_P(S)$. Now let $\bar C = C \backslash (C \cap S)$ and let $\bar S = S \backslash (C \cap S)$. The union of the lines joining a point of $\bar C$ and a point of $\bar S$ also contains an open set in $\PP^4$. Thus projection from a general point $P$ of $\PP^4$ has the property that $\pi_P(C) \not \subset \pi_P(S)$ and $| \pi_P(C) \cap \pi_P(S)| > bd = (\deg C)(\deg S)$,  a contradiction.
\vspace{.1in}

\item Assume that $C$ is a plane curve of degree $b \geq 2$ and assume $|C \cap S| = bd$.    Let $H$ be a general hyperplane. Let $P$ be a general point lying on the plane of $C$ and let $\pi_P$ be the projection from $P$ to $H$. Now   $\pi_P(C)$ is a line and the intersection points of $C$ and $S$ are mapped to distinct points in $H$. If the line $\pi_P(C)$ does not lie on $\pi_P(S)$ then  this line meets the surface in too many points: $bd = | C \cap S| \leq |\pi_P(C)\cap \pi_P(S)|=  d < bd$, a contradiction since $b \geq 2$. 

So assume that $\pi_P(C)$ lies on $\pi_P(S)$. This means that if $\Lambda$ is the plane of $C$ then $\Lambda \cap S$ is a curve (since $P$ is a general point of $\Lambda$). This curve is not necessarily equidimensional. Say it is the union of a plane curve $E$ of degree $e$ and a finite set  of  points.  If $Q$ is a general point of $S$ and $H$ is the hyperplane spanned by $Q$ and $\Lambda$, then $H$ meets $S$ in an equidimensional curve of degree $d$ having $E$ as a proper subcurve. So say $H \cap S = E \cup Y$ (as sets), where $Y$ is a curve of degree at most $d-e$ with no component on $\Lambda$. Thus $\Lambda \cap S$ consists at most of the curve $E$ plus $d-e$ points.  Hence $|C \cap S| \leq be + (d-e) < be + b(d-e) = bd$ and we have our contradiction.

\vspace{.1in}

\item 

Finally, assume that $C$ is a line, so $C$ is a $d$-secant line of $S$. Let $P$ be a general point of $S$ and $\pi_P$ the projection from $P$ to a general hyperplane. 
Then $S' := \overline{\pi_P(S)}$ (the  Zariski closure)  is a surface of degree $d-1$ inside the hyperplane. If the line $\pi_P(C)$ does not lie on $S'$ then we have a contradiction since it meets $S'$ in $d$ points, while a line meets a surface of degree $d-1$ in at most $d-1$ points.

 So assume that $\pi_P(C)$ lies on $\overline{\pi_P(S)}$ for all $P \in S$. This means that for every $Q \in C$, the line joining $Q$ and $P$ meets $S$ in at least one point other than $P$. Then the plane spanned by $C$ and $P$ meets $S$ in a curve. 
 
 Let $H$ be a general hyperplane containing $C$. Then $H \cap S$ is a curve of degree $d$, and $C$ is not a component of this curve. 
Furthermore, by Bertini and the irreducibility of $S$, $H \cap S$ is an irreducible curve of degree $d$. Then $C$ is a $d$-secant line to this curve. This is impossible unless $H \cap S$ is a plane curve of degree $d$. 

So finally assume that $H \cap S$ is a plane curve of degree $d$. Let $\Lambda$ be the plane of this plane curve. Consider the  pencil of hyperplanes containing $\Lambda$. Every hyperplane in this pencil meets $S$ in precisely this plane curve. But if $Q$ is a general point of $S$ and $H'$ is the element of the pencil containing $Q$ then $H' \cap S$ is more than the plane curve. Contradiction. (This last paragraph of the  proof is adapted from \cite[Proposition 18.10]{harris}.)
\end{itemize}

The last three bullet points show that $S$ cannot be nondegenerate, and then the first two show the desired conclusion.
\end{proof}

\begin{remark}
    In some examples above we find trivial geprofi sets Z which are degenerate, and for which at least one between the surface $S$ and the curve $C$ is nondegenerate. We notice that, in this case, $Z$ is also the complete intersection of a curve $C'$ and a surface $S'$ in a hyperplane $H$ containing $Z$. Indeed, one can take $S'$ as a generic projection of $S$ to $H$, and $C'$ as a generic projection of $C$ to $H$.
\end{remark}

Thus, we focus on the case of nondegenerate trivial geprofi sets, for which we provide a geometrical description. Namely, in the next result we show that for $Z$ to be a type 1 $(b,d)$-geprofi set, where by definition $Z$ is an intersection
of a curve $C$ of degree $b$ and a surface $S$ of degree $d$ meeting in $bd$ points, the curve $C$ is a union of plane curves and $S$ is a union of planes. 

\begin{theorem} \label{classi}
Let $C, S \subset \PP^4$ be such that $\deg S = d \geq 2, \deg C = b \geq 3$ and  $|C \cap S| = bd$, where $C \cap S$ is the set-theoretic intersection. Assume that $C\cap S$ is nondegenerate. Then each component of $S$ is a plane, i.e., $S=\pi_1\cup \dots\cup \pi_d$ is a nondegenerate union of planes. 
Moreover  $C$ is the union of (possibly reducible) plane curves $B_1,\dots, B_p$ of degree $>1$ and lines $L_1,\dots, L_q$. Each $L_j$ meets each plane $\pi_i$ in a point $P_{ij}$. If $\Pi_i$ is the plane spanned by $B_i$, then each $\Pi_j$ meets each plane $\pi_i$ in a line $\ell_{ij}$.

Furthermore, examples exist for each $d, p, q$, and for $1 \leq i \leq p$, the degree of $B_i$ can be chosen arbitrarily ($\geq 2$). Finally, for $d \geq 3$ and $p \geq 3$, this is essentially  the only  example.
\end{theorem}

\begin{proof}
It is arithmetically obvious that any component $S_i$ of $S$ meets any component $C_j$ of $C$ in $\deg(S_i)\deg(C_j)$ points. Thus, by Proposition \ref{triv irred}, all components $S_i$ and $C_j$ are degenerate.

Assume that $S$ has a component $S_0$ which is not a plane. Then $S_0$ spans a hyperplane $H_0$. By Proposition \ref{triv irred}, for all components $C_i$ of $C$ the union $S_0\cup C_i$ lies in a hyperplane, which must be $H_0$. Thus $C\subset H_0$ and $C\cap S\subset H_0$, a contradiction. Thus $S$ is a nondegenerate union of planes.

Similarly, if $C$ has a component $B_0$ which spans a hyperplane $H_0$, then, by Proposition~\ref{triv irred}, for all components $S_i$ of $S$ the union $B_0\cup S_i$ lies in a hyperplane, which must be $H_0$. Thus $S\subset H_0$ and $C\cap S\subset H_0$, a contradiction.

For arithmetical reasons, every $L_j$ must meet every $\pi_i$. Similarly, by Proposition \ref{triv irred}, every union $B_j\cup \pi_i$ is degenerate, and this forces the union $\Pi_j\cup \pi_i$ to be degenerate.

Next we address the  question  of existence. Suppose we have an example satisfying the following conditions:

Let $\pi_1,\dots,\pi_d$ be planes in $\PP^4$ whose union $S$ is nondegenerate. Let $\Pi_1,\dots,\Pi_p$ be planes and $L_1,\dots,L_q$ be lines in $\PP^4$. Assume that, for each $i,j$, the intersection $\pi_i\cap \Pi_j$ is a line $\ell_{ij}$, and the intersection $\pi_i\cap L_j$ is a point $P_{ij}$. Assume that all points $P_{ij}$ are distinct, and assume that these points, together with the planes $\Pi_1,\dots,\Pi_p$, span $\PP^4$.

    For $i=1,\dots,p$ fix a general plane curve $B_i\subset\Pi_i$, of degree $b_i>1$, and define $C=B_1\cup\dots\cup B_p\cup L_1\cup\dots\cup L_q$. 

    Then $C$ is a curve of degree $b=q+\sum b_i$, which meets $S$ in a nondegenerate set $Z$ of $bd$ points. Clearly $Z$ is a trivial geprofi set. Namely, the generality of the curve $B_i$ guarantees that it meets each line $\ell_{ij}$ in $b_i$ points, all distinct and none of them matching with some $P_{hk}$.

    Notice that we could take as $B_i$ any (even reducible) curve in $\Pi_i$, provided that the intersection $C\cap S$ has cardinality $bd$ and spans $\PP^4$ (so, in particular,  each $B_j$ meets the union of the lines $\ell_{1j},\dots,\ell_{dj}$ in $b_id$ distinct points).

    Now we prove the existence. Fix $d \geq 2$, $p \geq 1$ and $q \geq 1$. In a hyperplane $\PP^3 = H \subset \PP^4$ fix a smooth quadric surface $\mathcal Q$. Let $\ell_1,\dots,\ell_d$ be lines in one ruling of $\mathcal Q$ and $m_1,\dots,m_p$ be lines in the other ruling. Fix a point $P \in \PP^4$ not on $H$. For $1 \leq i \leq d$ let $\pi_i$ be the plane spanned by $\ell_i$ and $P$, and for $1 \leq j \leq p$ let $\Pi_j$ be the plane spanned by $m_j$ and $P$. It is clear that each $\pi_i$ meets each $\Pi_j$ in a line. Furthermore, any two $\pi_i$ meet in a point, as do any two $\Pi_j$.  

Fix any  $j$. On the plane $\Pi_j$ fix a general curve $B_j$ of degree $b_j >1$. Note that since $\Pi_j$ contains the line spanned by $P$ and the intersection of  $\ell_i$ and $m_j$, $B_j$ meets this line in $b_j$ distinct points, so $B_j$ meets $\pi_i$ in $b_j$ points for each $i$.

For the $L_k$, $1 \leq k \leq q$, observe that the cone over $\mathcal Q$ with vertex $P$ has a general hyperplane section that  is a smooth quadric surface in the hyperplane. Choosing a line $L$ in the  appropriate ruling, $L$ meets each $\pi_i$ in a point. Then taking $p$ general hyperplane sections and choosing $L_i$ in the corresponding quadric surface, we produce the $L_i$. 

Let $S = \pi_1 \cup \dots \cup \pi_d$ and let $C = B_1 \cup \dots \cup B_p \cup L_1  \cup \dots \cup L_q$ as above. As long as $d \geq 2$, $S$ is nondegenerate. If $p \geq 2$ and $q \geq 0$ then $C$ is also nondegenerate. If $p \geq 1$ and $q \geq 1$ then again $C$ is nondegenerate. These choices also give $C \cap S$ nondegenerate.

Finally, the  fact that this example is essentially unique follows from Proposition  \ref{one point} below.
\end{proof}

\begin{remark}  
    We stress that one can use the construction of Example \ref{triv ex 1}, with $S=$ union of planes and $C=$ one plane curve of degree $b-1$ plus one general line, to prove that for all $d\geq 2$ and $ b\geq 3$ there exists a nondegenerate trivial $(b,d)$-geprofi set.
\end{remark}

It follows from Theorem \ref{classi} and its proof that a complete classification of trivial geprofi sets in $\PP^4$ would be based on the classification of configurations of lines $L_1,\dots,L_q$ and planes $\pi_1,\dots,\pi_d$, $\Pi_1,\dots,\Pi_p$, such that each $\pi_i$ meets each $L_j$ in a point and each $\Pi_k$ in a line. 

The classification of such configurations is far from obvious.

The following proposition goes in the direction of having a complete classification.

\begin{proposition} \label{one point}
      Let $\pi_1, \ldots, \pi_d$ be planes in $\mathbb P^4$ such that the union of any two of them is nondegenerate, $d\ge 3$. 
    Let $\Pi_1,\Pi_2$ be distinct planes such that    $\pi_i\cap \Pi_j=\ell_{ij}$ for all $i,j$. Assume that the lines $\ell_{ij}$ are distinct.  
    Then, the planes $\pi_1, \ldots, \pi_d$ have a point in common.
\end{proposition}
\begin{proof}
Since any two of the planes $\pi_1,\ldots,\pi_d$ span $\mathbb P^4$ we can introduce the points $Q_{ij}=\pi_i\cap \pi_j.$ 

Assume by contradiction that there is a plane say $\pi_1$ which meets the other planes in (at least two) different points, say  $Q_{12}\neq Q_{13}$. 
We collect the following relevant facts
\begin{itemize}
    \item[(1)] The lines $\ell_{11},\ldots, \ell_{d1}$ are all contained in the plane $\Pi_1$, so they (since are distinct lines) meet pairwise in a point.
\item[(2)] since  $\ell_{i1}\cap \ell_{j1}\subseteq \pi_i\cap \pi_j$, the points $\ell_{i1}\cap \ell_{j1}$  are precisely the $Q_{ij}$s.
\item[(3)] The plane $\Pi_1$ contains all the points $Q_{ij}$. The same is true for $\Pi_2$.
\item[(4)] Since $\Pi_1\neq \Pi_2$ then the points $Q_{ij}$ must lie on a line $V$, which is the span of $Q_{12}$ and $Q_{13}$. Therefore, $V\subseteq \pi_1.$
\end{itemize}

Thus $\pi_1$ contains the point $Q_{23}$ and hence either $\pi_1\cap \pi_2=V$ or $\pi_1\cap \pi_3 =V$ which contradicts the hypotheses. 
\end{proof}


\section{First examples of  geprofi sets of type 2} \label{sec: 10p}
In this section we initiate the study of nontrivial geprofi sets.
We begin with a lemma that will also be useful later in this paper.

\begin{lemma} \label{l. 2b points on a quadric}
Let $Y$ be a set of $2b\geq 4$ points in $\PP^3$, contained in a given quadric $\mathcal Q$. Assume no point of $Y$ is a vertex for $\mathcal Q$ and one of the following conditions holds:
\begin{itemize}
    \item[(i)] $\mathcal Q=\Pi_1\cup \Pi_2$ is reducible and the points of $Y$ are equidistributed on the two planes; note that the assumption that no point of $Y$ is a vertex for $\mathcal Q$ means that the points of $Y$ avoid the intersection  line of the planes; 

     \item[(ii)] no $b+1$ points of $Y$ are on a line.
\end{itemize}
Then there exists a partition of $Y$ in sets of cardinality 2 such that all the lines $\ell_{ij} = \langle P_i,P_j\rangle$, which join two points of the same subset of the partition,  do not lie on $\mathcal Q$.  Hence  if $Y$ arises as the image of the general projection of a set $Z$ of points in $\PP^4$ then $Z$ is a $(b,2)$-geprofi set.
\end{lemma}

\begin{proof}
When $b = 2, 3, 4$ there is an infinite family of such quadrics (i.e., quadrics containing $Y$); when $b \geq 5$ it is a nontrivial assumption that $Y$ lies on a quadric surface.

For any value of $b$, if the assumption of (i) holds then the result is trivial  since the lines spanned by a point in $\Pi_1$ and a point in $\Pi_2$ are not in $\mathcal Q$. So from now on we assume (ii). Our proof will be by induction  on $b$.

First assume that $b=2$. Denote by $P_1,P_2,P_3,P_4$ the four points and by $\ell_{ij}$ the line joining $P_i$ and $P_j$.   By assumption there are no three points on a line, hence there are 6  distinct lines $\ell_{ij}$ and at most 4 of them are contained in $\mathcal Q$. If the two lines not in $\mathcal Q$ cover $Y$ we are done. 
Say $\ell_{12}\not \subseteq \mathcal Q$ and $\ell_{13}\not \subseteq \mathcal Q$ and $\ell_{34},\ell_{24}\subseteq \mathcal Q$ (otherwise we are done again). Then   $\ell_{14}\not \subseteq \mathcal Q$ (because there are already two lines  in $\mathcal Q$ containing $P_4$, and $P_4$ is not a vertex of $\mathcal Q$) and $\ell_{23}\not \subseteq \mathcal Q$ (because, otherwise, $\mathcal Q\cap \langle P_2,P_3,P_4\rangle$ would be a cubic curve) so we are done.

Now assume $b>2$. 
Consider  a line $L_1$ in $\mathbb P^3$ which intersects $Y$ in the maximal number of points and a line $L_2$ which intersect $Y\setminus {L_1}$ in the maximal number of points.  Note that there is at most one set of $b$ collinear points in  $Y\setminus {L_1}$ since $|Y| = 2b$.

Fix any set $Y'$ of 4 points, two of them in $L_1$ and the other two in $L_2\setminus L_1$.

For one of the points in $Y'$ there are three lines joining it with the other three points. But at most two of the lines are in $\mathcal Q$. Hence there are $P_1,P_2\in Y'$ such that  the line  $\langle P_1,P_2 \rangle\not \subseteq \mathcal Q$. Note that if $|L_i\cap Y|>2$ then  $L_i$ lies on $\mathcal Q$ so $P_1,P_2$ cannot both lie on $L_i$. Therefore we consider the set  $Y\setminus \{P_1,P_2\}$, which has no $b$ points on a line, and  use induction.
\end{proof}

\begin{remark}
Using Lemma \ref{l. 2b points on a quadric} it is easy to check that a set $Z$ of 8 general points in $\PP^4$ is a nontrivial $(4,2)$-geprofi set (it is nontrivial since clearly $Z$ is not contained in two planes and the irreducible quadric surfaces of $\mathbb P^4$ are degenerate).
\end{remark}
In Theorem \ref{2b LGP class} we classify $(b,2)$-geprofi sets under the extra assumption that they are in linear general position (LGP). Without this latter assumption, it is more delicate to describe geprofi sets. One of the first interesting cases is that of 10 points and the question of whether they are $(5,2)$-geprofi.  
The general projection of 10 general points in $\mathbb P^4$ is not contained in a quadric surface. Indeed in this case, the 2-Weddle locus (i.e., the locus of points in $\PP^4$ such that the projection of the $10$ points from them is contained in a quadric surface in $\PP^3$) is a threefold of degree $5$; see \cite[Example 2.8]{POLITUS1} and Remark \ref{classical weddle} below for the definition and further details.

In the following example we give a nontrivial (5,2)-geprofi set.

\begin{example}\label{ex. 10 points}
Let $Z=\{P_1, \ldots,  P_{10}\}$ be the set containing the following points
    $$
\begin{array}{lr}
    P_1=(1:0:0:0:0), & P_{10}=(0:1:1:1:1), \\
    P_2=(0:1:0:0:0), & P_{9}=(1:0:1:1:1), \\
    P_3=(0:0:1:0:0), & P_{8}=(1:1:0:1:1), \\
    P_4=(0:0:0:1:0), & P_{7}=(1:1:1:0:1), \\
    P_5=(0:0:0:0:1), & P_{6}=(1:1:1:1:0). \\
\end{array}$$
We note that lines $\ell_i=P_iP_{11-i}$ for $i=1,\ldots,5$ all intersect  in the point $O=(1:1:1:1:1)$.

We note that for a general point $Q=(a_0:a_1:a_2:a_3:a_4)\in\PP^4$ there exists a quadratic cone $C_Q$ with vertex at $Q$ containing all points in $Z$. One can check that 
$$C_Q:\; \sum_{0\leq i<j\leq 4}(-1)^{i+j}(a_j-a_i)(a_m-a_l)(a_m-a_k)(a_l-a_k)x_ix_j=0,$$
where for fixed $i<j$ the values for the letters $k,l,m$ are taken so that $k<l<m$ and $\left\{i,j,k,l,m \right\}=\left\{1,2,3,4,5\right\}.$

Hence, from Lemma \ref{l. 2b points on a quadric}, the set $Z$ is a $(5,2)$-geprofi set. 
This example is nontrivial since $Z$ spans $\mathbb P^4$ and it is not contained in two planes.

We can say a little bit more. We claim that the set $Z$  has $h$-vector $h_Z=(1,4,5)$. 
 To see this it suffices to show that for any index $1\leq i\leq 10$ there exists a quadric $\mathcal{Q}_i$ vanishing in all points of $Z$ with the exception of $P_i$. By symmetry it is enough to do this for $P_1$ and $P_{10}$. For $P_1$ we can in fact show two planes containing the remaining $9$ points in $Z$:
$$\mathcal{Q}_1\colon\; x_0(x_1+x_2+x_3+x_4-3x_0)=0.$$
For $P_{10}$ it is easy to verify that the remaining $9$ points are also contained in the union of two hyperplanes, so that
\[
\mathcal{Q}_{10}\colon\; x_1(x_0+x_2+x_3+x_4-3x_1)=0
\]
does the job. It follows that $C_Q$ is an unexpected cone  for $Z$.

\end{example}

The next   proposition and remark generalize   Example \ref{ex. 10 points}.
\begin{proposition}\label{p. 10 points}
    Let $\ell_1, \ldots, \ell_5\subseteq \mathbb P^4$ be five lines which are concurrent in a point $O$ and such that the union of any four of them spans $\mathbb P^4$.  Let $Z$ be a set of 10 points equidistribuited on the lines, such that $O\notin Z$. Then $Z$ is a nontrivial $(5,2)$-geprofi set. 
\end{proposition}
\begin{proof}
    Let $Q$ be a general point in $\mathbb P^4$. The projection from $Q$ of $\ell_1\cup \cdots\cup \ell_5$ is contained in a quadric cone $\mathcal Q$ with vertex at $\pi_Q(O)$. Such cone is unique since among the five lines $\ell_i$ no subset of four is contained in a plane. 
    Thus, in particular, $\pi_Q(Z)$ is contained in a quadric cone and, from Lemma \ref{l. 2b points on a quadric}, $Z$ is a $(5,2)$-geprofi set. Since $\ell_1\cup \cdots\cup \ell_5$ spans $\mathbb P^4$ so does $Z$. 

    To show that $Z$ is  nontrivial it is enough to show that no set of five points of $Z$ are contained in a plane.  Let $P_i,Q_i$ the points of $Z$ contained in $\ell_i$.  
    Since four lines span $\mathbb P^4$ then any three of them span $\mathbb P^3$, hence there cannot be 5 coplanar points on three lines. 
    Say $P_1, P_2,P_3,P_4$ are on a same plane $\pi_1$,  then $\pi_1$ together with $O$ span a hyperplane $H$ and the lines $\ell_1, \ldots, \ell_4$ would be on $H.$
\end{proof}

\begin{remark}\label{unexpected cones 10 points}

    There are two possibilities  for the $h$-vector of the sets $Z$ constructed in Proposition \ref{p. 10 points}. We have either $h_Z=(1,4,5)$ or $h_Z=(1,4,4,1)$ depending on whether the 10 points impose independent conditions on quadrics or not.  To see this, note first that by construction the 10 points are nondegenerate in $\PP^4$, hence the degree 1 entry must be 4. Let $C = \ell_1 \cup \dots \cup \ell_5$. By construction, $C$ is a cone over a set of 5 points in $\PP^3$ in LGP. Such a set of 5 points is arithmetically Gorenstein, so the same is true of $C$. The first difference of the Hilbert function of $C$ is $(1,4,5,5,\dots)$. If $Z$ is cut out on $C$ by a quadric hypersurface then the $h$-vector of $Z$ is $(1,4,5-1, 5-4, 5-5, \dots) = (1,4,4,1)$. Otherwise the $h$-vector of $Z$ is $(1,4,5)$ as claimed. For the case $(1,4,4,1)$, since $C$ is a cone over an arithmetically Gorenstein set of points in $\PP^3$, $C$ is also arithmetically Gorenstein, and hence so is $Z$. Note that in this latter case, $Z$ is a full intersection of a quadric hypersurface and the curve  $C = \ell_1\cup \cdots\cup \ell_5$.

We now focus on the case  where the $h$-vector of $Z$ is $(1,4,5)$ and on some interesting properties that it possesses. 

\begin{enumerate}

\item \label{quad cont Z} All the quadrics vanishing at $Z$ contain the curve $ C = \ell_1 \cup \ldots \cup \ell_5$. Indeed, the $h$-vector shows that the vector space of quadrics containing $C$ and the vector space of quadrics containing $Z$ coincide. 

\item Let $W$ be any subset of 9 of the points of $Z$. Notice that the $h$-vector of $W$ is $(1,4,4)$. In contrast to (\ref{quad cont Z}), a general quadric hypersurface containing $W$ does not contain $C$, and instead cuts out a set $Z'$ of 10 points on $C$. Setting $P_W$ to be the tenth point, we have $Z' = W \cup P_W$.  Since $C$ is arithmetically Gorenstein, so is $Z'$, and so $P_W$ is the point residual to $W$ in this Gorenstein link. As above, the $h$-vector of $Z'$ is $(1,4,4,1)$.  These $h$-vectors then imply
\[
[I_Z]_2 \subsetneq  [I_W]_2  =  [I_{Z'}]_2.
\]
In particular, $P_W$ depends only on $W$ and not on the quadric chosen, and furthermore $P_W \notin Z$. Hence each subset $W \subset Z$ of 9 points has associated to it a tenth point that is not a point of $Z$.

\item Any singular quadric in $[I_W]_2$ must in fact contain all of $C$. Indeed, note first that a general projection of $W$ lies on a unique quadric surface in $\PP^3$, as does a general projection of $C$, as was discussed in the proof of Proposition \ref{p. 10 points}. Thus these quadrics are the same. Therefore, a general singular quadric containing $W$ must contain all of $Z$. The conclusion then comes by semicontinuity.

\item The existence of a quadric cone $C_Q$ containing $Z$ and having vertex at a general point $Q$ is unexpected in the sense of \cite{CHMN} (this is not true for the case $h_Z=(1,4,4,1)$). Indeed $\dim[I_Z]_2=5$ and a double point imposes at most five conditions on the linear system of quadric forms.    
See   \cite[Definition 2.5]{CM}, \cite{HMNT, HMTG} for more on unexpected cones.
\end{enumerate}
\end{remark}

 In   Remark \ref{conn to lef} we show a relation between the results of this section and the following question about the Weak Lefschetz Property, via Macaulay duality. (This is also related to Remark \ref{classical weddle}.)
\begin{question}\label{q. 10 pts WLP}
Let $L_1,\dots,L_{10}$ be linear forms in $R = k[x_0,\dots,x_4]$. Does the quotient algebra $R/(L_1^2,\dots,L_{10}^2)$ have the Weak Lefschetz Property (WLP)? 
\end{question}
\noindent Following the ideas in \cite{POLITUS1}, we show that the answer is deeply connected to the  geometry of ten points in $\mathbb P^4$ studied above. 

\begin{remark} \label{conn to lef}
Let $L_1,\dots,L_{10}$ be linear forms in $R = k[x_0,\dots,x_4]$ and 
let $P_1,\dots,P_{10}$ be their dual points in $\PP^4$. Set $\wp_i = I(P_i)$. For the purposes of this remark, we will assume that $P_1,\dots,P_{10}$ impose independent conditions on quadrics (i.e., the $h$-vector is $(1,4,5)$), so $\dim [\wp_1 \cap \dots \cap \wp_{10}]_2 = 15 - 10 = 5$.

By Macaulay duality we have
\[
5 = \dim [\wp_1 \cap \dots \cap \wp_{10}]_2 = \dim [R/(L_1^2,\dots,L_{10}^2)]_2.
\]
Of course we also have $\dim [R/(L_1^2,\dots,L_{10}^2)]_1 = \dim [R]_1 = 5$.
Let $I = (L_1^2,\dots,L_{10}^2)$. For a  linear form $L$ we have the homomorphism
\[
\times L : [R/I]_1 \rightarrow [R/I]_2.
\]
To answer Question \ref{q. 10 pts WLP}, we want to know whether this fails to be an isomorphism for all $L$. That is, we want to know whether the cokernel of this homomorphism is nonzero, i.e.,
\[
\dim [R/(I,L)]_2 > 0.
\]
Let $P$ be the point dual to $L$, whose ideal is $\wp$. Invoking Macaulay  duality again, we have
\[
\dim [R/(I,L)]_2 = \dim [R/(L_1^2,\dots,L_{10}^2, L)]_2 = \dim[\wp_1 \cap \dots \cap \wp_{10} \cap \wp^2]_2.
\]
So $R/I$ fails to be an isomorphism  from degree 1 to degree 2 if and only if $Z = \{ P_1,\dots,P_{10} \}$ has the  property that the general projection lies on a quadric surface in $\PP^3$. 

The conclusion from this section  is that if $L_1,\dots,L_{10}$ are sufficiently general then for a general linear form $L$, the multiplication $\times L : [R/I]_1 \rightarrow [R/I]_2$ is an isomorphism  (because the general projection of $\{ P_1,\dots,P_{10} \}$ does not lie on a quadric surface in $\PP^3$). But if we move these points slightly, so that two points lie on each of five sufficiently general concurrent lines (maintaining the property that they impose independent conditions on quadrics, i.e., they have $h$-vector $(1,4,5)$), then the corresponding isomorphism fails.
Notice that the special points do not have LGP, but they are close: at least they have the property that no three lie on a line and at most four lie on a plane.

Corollary \ref{c. LGP + 1441} will show that if $Z$ is in LGP then it  has $h$-vector $(1,4,4,1)$ if and only if it is a $(5,2)$-geprofi set. Thus, the analysis above  implies that the map  $\times L$ is an isomorphism if and only if $Z$ has $h$-vector $(1,4,4,1)$.
\end{remark}


\section{A useful construction for nontrivial geprofi sets} \label{further ex}

 Proposition \ref{p. 10 points} can be generalized and also extended to a larger class of nontrivial geprofi sets not in LGP. The idea behind the examples constructed in this section is contained in the following trivial lemma.
\begin{lemma}\label{l.geprofi sets not LGP}
Let $Z\subset \mathbb P^4$  be a set of $bd$ points such that 
\begin{itemize}
    \item there exist $b$ lines in $\mathbb P^4$, $\ell_1, \ldots, \ell_b$,  each containing $d$ points of $Z;$
    \item for  a general point $P\in \mathbb P^4$, there exists a point $Q\in \mathbb P^3$ such that  $\pi_Q(\pi_P(Z))$ lies on a curve $\mathcal C$ of degree $d.$
\end{itemize}
Then,
\begin{itemize}
    \item[(i)] $\pi_P(Z)$ lies on a surface $\mathcal S$ of degree $d$ (the cone with vertex $Q$ over $\mathcal C$);
    \item[(ii)]if none of the lines $\pi_P(\ell_1),\dots, \pi_P(\ell_b)$ are  contained in $\mathcal S$  then $Z$ is $(b,d)$-geprofi. 
\end{itemize}
\end{lemma}

Despite Lemma \ref{l.geprofi sets not LGP} being so immediate, it gives a hint to construct a special class of $(b,d)$-geprofi sets that are nontrivial and not in LGP, with $b-2\ge d$. 

\begin{proposition}\label{p.geprofi sets not LGP}
Let $H$ be a hyperplane in $\mathbb P^4$ and let $X=\{P_{ij}\}_{1\le i,j\le d}$ be a $(d,d)$-grid on a quadric surface contained in $H$, where $d\ge 3$. Call $L_1,\ldots, L_d$ and $M_1,\ldots, M_d$  the rulings defining $X$.
Let $P_0$ be a general point in $\mathbb P^4$.
In the plane spanned by $L_1$ and $P_0$ consider a general line $L$ and let $X_1$ be the set of $d$ collinear points obtained as the intersection of $L$ with the lines $P_{1j}P_0$.
In the plane spanned by $M_1$ and $P_0$ consider a general line $L'$ and let $X_2$ be the set of $d$ collinear points obtained as the intersection of $L'$ with the lines $P_{j1}P_0$. (See the figure below.)

Set $Z=X\cup X_1\cup X_2$.
Then 
\begin{enumerate}
    \item $Z$ is $(d+2,d)$-geprofi;
    \item $Z$ is not a trivial geprofi set. 
\end{enumerate}   
\end{proposition}
    \begin{proof}
    
(1) Let $P$ be a general point in $\mathbb P^4$ and let $Q=\pi_P(P_0)$.  Note that $\pi_Q(\pi_P(Z)) = \pi_Q(\pi_P(X))$, is a complete intersection of $d^2$ points in the plane. Then,  $\pi_Q(\pi_P(Z))$ is the base locus of a pencil of curves of degree $d$ in $\mathbb P^2$. The general element in such a linear system is irreducible (Bertini).
 Then, the set $\pi_P(Z)$ is contained in an irreducible cone of degree $d$ with vertex at $Q$. The lines $\pi_P(L_1),\ldots,\pi_P(L_d),\pi_P(L),\pi_P(L')$  do not contain the point $Q$. Hence, by Lemma \ref{l.geprofi sets not LGP}, $Z$ is $(d+2,d)$-geprofi. 

(2) Assume by contradiction that $Z$ is a trivial geprofi set given by the intersection in $\mathbb P^4$ of a curve $C$ of degree $d+2$ and a surface $S$ of degree $d$.  From Theorem \ref{classi} the surface $S$ is a union of planes, $S=\pi_1\cup \cdots \cup \pi_d$. 
The set $Z$ is contained in at least two curves of degree $d+2$, namely $L_1\cup \cdots\cup L_d\cup L\cup L'$ and $M_1\cup \cdots\cup M_d\cup L\cup L'$.
We have the following cases.

\begin{itemize}
    \item If both $L$ and $L'$ are not contained in any of the planes $\pi_i$ then  each plane $\pi_i$ contains exactly one point of $X_1$ and one point of $X_2$. Therefore, because of the generality of $L$ and $L'$, none of the planes $\pi_i$ contains any of the rulings $L_j,M_k$ defining $X$. Thus each $\pi_i$ contains at most 4 points of $Z$ (one from $X_1$, one from $X_2$ and two from $X$) and this contradicts $d\ge3$.
    
    \item Assume $L\subseteq \pi_1$ and $L'$ is not contained in any of the planes $\pi_i$. 
    From the generality of $L'$, then, the planes $\pi_i$ must be the span of the rulings $L_i$ and the point $P_0$.   The curve $C$ is either a union of $d+2$ lines or else a union of lines and plane curves of smaller degree, but the collinearities of the points of $Z$ force it to be a union of lines. Then $L$ must be a component; but such a curve does not meet $S$ transversally.
    
\item Finally, if we assume $L\subseteq \pi_1$ and $L'\subseteq \pi_2$, by the generality of $L$ and $L'$, $\pi_1\cup \pi_2$ contains at most $2d-1$ points of $X$. This case happen if, say, $L_1\subseteq \pi_1$ and $M_1\subseteq \pi_2$. The residual $(d-1)^2$ points are on a $(d-1,d-1)$-grid which cannot be contained in $d-2$ planes.  Any other possibility contains even fewer points on $\pi_1 \cup \pi_2$. \qedhere
\end{itemize}
\end{proof}

\definecolor{uququq}{rgb}{0.35,0.35,0.35}

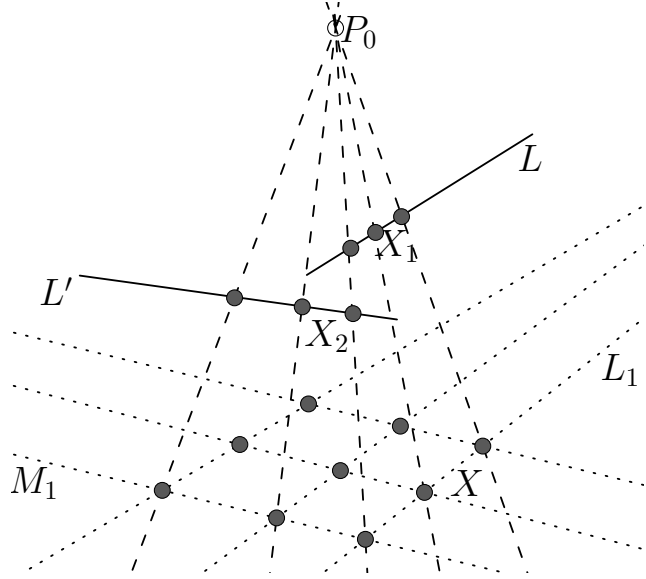
\begin{figure}[ht] \label{FigureForProp}
\centering
\begin{tikzpicture}[scale=0.3,line cap=round,line join=round, x=1cm,y=1cm]
\clip(-16,-6) rectangle (12,19.3);
\draw [line width=0.7pt,loosely dotted,domain=-18.02065158759578:25.10898725734244,color=black] plot(\x,{(--3.3023545041968205-0.9855823628699303*\x)/4.0708836727236255});
\draw [line width=0.7pt,loosely dotted,domain=-18.02065158759578:25.10898725734244,color=black] plot(\x,{(-8.145877742049434-1.1569879911951357*\x)/4.456546336455338});
\draw [line width=0.7pt,loosely dotted,domain=-18.02065158759578:25.10898725734244,color=black] plot(\x,{(-9.628431951673537-1.7997590974146553*\x)/-3.042449902772394});
\draw [line width=0.7pt,loosely dotted,domain=-18.02065158759578:25.10898725734244,color=black] plot(\x,{(--1.1069885584456587-1.9711647257398608*\x)/-2.656787239040682});
\draw [line width=0.7pt,loosely dotted,domain=-18.02065158759578:25.10898725734244,color=black] plot(\x,{(-23.212040546337995-1.2220449824282644*\x)/5.065505819311399});
\draw [line width=0.7pt,loosely dotted,domain=-18.02065158759578:25.10898725734244,color=black] plot(\x,{(--11.034462147360832-2.0539699926047126*\x)/-2.5841789586170165});
\draw [line width=0.7pt,dash pattern=on 4pt off 6pt, domain=-18.02065158759578:25.10898725734244,color=black] plot(\x,{(-172.76443813940284-20.502065403776477*\x)/-7.679211576601889});
\draw [line width=0.7pt,dash pattern=on 4pt off 6pt,domain=-18.02065158759578:25.10898725734244,color=black] plot(\x,{(-82.73880509881677-21.724110386204742*\x)/-2.6137057572904903});
\draw [line width=0.7pt,dash pattern=on 4pt off 6pt,domain=-18.02065158759578:25.10898725734244,color=black] plot(\x,{(-12.796535368689764-22.673535532094164*\x)/1.3217619797048252});
\draw [line width=0.7pt,dash pattern=on 4pt off 6pt,domain=-18.02065158759578:25.10898725734244,color=black] plot(\x,{(--38.12582275913831-20.593029803728765*\x)/3.9393265721193593});
\draw [line width=0.7pt,dash pattern=on 4pt off 6pt,domain=-18.02065158759578:25.10898725734244,color=black] plot(\x,{(--88.39869349495001-18.53905981112405*\x)/6.523505530736376});
\draw [line width=0.7pt,domain=-12.96169577457285:1.1,color=black] plot(\x,{(--15.310270873675325-0.39336501456972783*\x)/2.8397866322740626});
\draw [line width=0.7pt,domain=-2.9:7.1,color=black] plot(\x,{(--22.55764460424735--1.5606063634635472*\x)/2.500271620644382});
\begin{scriptsize}
\draw [ color=black, fill=uququq] (-2.8230923096666363,1.4946986901463075) circle (10pt);
\draw [color=black, fill=uququq] (1.2477913630569892,0.5091163272763772) circle (10pt);
\draw [color=black, fill=uququq] (-5.86554221243903,-0.3050604072683479) circle (10pt);
\draw [color=black, fill=uququq] (-1.4089958759836927,-1.4620483984634836) circle (10pt);
\draw [color=black, fill=uququq] (-9.302464487992086,-2.338169359003957) circle (10pt);
\draw [color=black, fill=uququq] (-4.236958668680687,-3.5602143414322214) circle (10pt);
\draw [color=black, fill=uququq] (4.900252619346179,-0.3751637663515317) circle (10pt);
\draw[color=black] (4.2,-2.0) node {{\large$X$}};
\draw [color=black, fill=uququq] (2.3160736607291623,-2.4291337589562443) circle (10pt);
\draw [color=black, fill=uququq] (-0.30149093168537183,-4.509639487321644) circle (10pt);
\draw [color=black] (-1.623252911390197,18.16389604477252) circle (10pt);
\draw[color=black] (-0.6,18) node {{\large $P_0$}};
\draw [color=black, fill=uququq] (-6.1,6.2) circle (10pt);
\draw [color=black, fill=uququq] (-3.1,5.8) circle (10pt);
\draw[color=black] (-14.0,6.5) node {{\large$L'$}};
\draw[color=black] (-15,-2) node {{\large$M_1$}};
\draw[color=black] (11,3) node {{\large$L_1$}};
\draw [color=black, fill=uququq] (-0.95,8.4) circle (10pt);
\draw [color=black, fill=uququq] (1.3,9.8) circle (10pt);
\draw[color=black] (1.1,8.5) node {{\large$X_1$}};
\draw[color=black] (7,12.4) node {{\large$L$}};
\draw [color=black, fill=uququq] (-0.85,5.5) circle (10pt);
\draw[color=black] (-2,4.5) node {{\large$X_2$}};
\draw [color=black, fill=uququq] (0.15,9.1) circle (10pt);
\end{scriptsize}
\end{tikzpicture}
\caption{An example of $(5,3)$-geprofi as in Proposition \ref{p.geprofi sets not LGP}.}
\end{figure}

The case $d=2$ is excluded in the assumption of Proposition \ref{p.geprofi sets not LGP}; however it is possible to construct, in a similar way, a class of $(b,2)$-geprofi sets, with $b\ge 5$.
\begin{proposition}\label{p.5 lines geprofi}
    Take a point $P_0\in \mathbb P^4$ and five lines $L_1,\ldots, L_5$ through $P_0$ not contained in a hyperplane.  Let $Z$ be a non degenerate set of $2b$ points on the five lines such that
   \begin{itemize}
        \item no three lines lie on a plane;
       \item each line contains at least one point of $Z$;
       \item each line contains at most $b$ points of $Z$.
   \end{itemize}
Then $Z$ is a nontrivial $(b,2)$-geprofi set. 
\end{proposition}
\begin{proof}
    For a general point $P\in \mathbb P^4$, the above first and second conditions ensures that $\pi_Q(\pi_P(Z))$ (where $Q=\pi_P(P_0)$) lies on a irreducible conic. Hence the quadric cone   containing $\pi_P(Z)$ with vertex at $Q$ cannot split in the union of  two planes.  Using the third condition we can find $b$ lines containing $Z$ without using $L_1,\ldots, L_5$.
\end{proof}

The method in this section gives an existence theorem for non trivial geprofi sets for some values of the parameters. In  Section \ref{LGPcurves} we will construct a class of geprofi sets in LGP which cover the missing cases of the following result, leading to a complete numerical classification in Theorem \ref{numerical classification}.

\begin{theorem} \label{t.extend} 
    There exists a nondegenerate nontrivial $(b,d)$-geprofi set for each $d\ge 2$ and $b\ge d+2$. 
\end{theorem} 

\begin{proof}We have the following cases.
\begin{itemize}
    \item A set of eight general points in $\mathbb P^4$ gives a nontrivial  $(4,2)$-geprofi set. 

\item For $d=2$ and $b\ge 5$ we use the construction in Proposition \ref{p.5 lines geprofi}.

\item For $d>2$ and $b=d+2$ we use the construction in Proposition \ref{p.geprofi sets not LGP}

\item The cases $d>2$ and $b>d+2$ can be obtained from a $(d+2,d)$-geprofi set as constructed in Proposition \ref{p.geprofi sets not LGP} by adding 
more sets of collinear points, for instance, in the  plane spanned by $P_0$ and $L_1$.\qedhere
\end{itemize} 
\end{proof}


\section{On $(b,2)$-geprofi sets in linear general position}
\label{sec: b,2 geprofi}

In \cite{POLITUS1} the authors conjectured that geproci sets in $\PP^3$ in linear general position (LGP) do not exist apart from sets of four points. In this section we begin studying to what extent geprofi sets of points in $\PP^4$ in LGP do exist. We will find examples of such sets which are not trivial  geprofi sets (indeed, Theorem \ref{classi} excludes sets in LGP from being trivial geprofi sets except for very few points). However, most of our examples do lie on curves of degree $b$ in~$\PP^4$.

\begin{proposition}\label{thequad}
Let $Z\subset\PP^4$ be a set of $2b$ points in LGP. Assume that a general projection $Z'$ of $Z$ to a hyperplane $L$ is contained in a quadric surface in $L$. Then $Z$ is a $(b,2)$-geprofi set.
\end{proposition}

Before we prove the proposition, we need the following observations.

\begin{remark} \label{gen proj 6 pts}
Let $Z$ be a set of $6$ points in $\PP^3$ in LGP and let $Z = A \cup B$ be a partition of $Z$ into two sets of three points. Let $\Lambda$ be the plane spanned by $A$ and let $P$ be a general point of $\Lambda$. Let  $\pi_P$ be the projection from $P$ to a general plane $\PP^2$. Then $\pi_P(A)$ is a set of three collinear points in $\PP^2$, while $\pi_P(B)$ is a set of three noncollinear points. Thus $\pi_P(Z)$ does not lie on a conic, which by semicontinuity means that projection from a general point of $\PP^3$ sends $Z$ to six points in the plane not lying on a conic. This was known already to Weddle in the 1800’s, and is the starting point to studying the so-called Weddle surface \cite{emch}, \cite{POLITUS1}.
\end{remark}

\begin{lemma}\label{nocon}
Let $Z$ be a set of $q\geq 7$ points in LGP in $\PP^4$. Then a general projection $Z'$ of $Z$ to $\PP^3$ is not contained in a quadric cone with vertex in a point of $Z'$. \end{lemma}

\begin{proof} It is enough to prove the statement for $q=7.$
Let $Q \in \PP^4$ be a general point. 
Assume by contradiction that there exists a point $P_1\in Z$ such that the  projection $Z' = \pi_Q(Z)$ to $\PP^3$ lies in a quadric cone with vertex at the projection $P'_1 = \pi_Q(P_1)$. Let $Y = Z \backslash \{P_1\}$, a set of 6 points, and note that $Y$ is in LGP. 

We will consider the projections $\pi_Q$ and $\pi_{P_1}$. Note that $P_1'$ has the property that $P_1' = \pi_{P_1}(Q)$ (both are the intersection of the line $\lambda = \overline{P_1Q}$ with the target $\PP^3$). 
We will also consider the projection $\pi_\lambda : \PP^4 \backslash \lambda \rightarrow \PP^2$.

Since the points of $Z$ (including $P_1$) are fixed, varying $Q$ allows us to conclude that $\pi_{P_1} (Q)$ is a general point of $\PP^3$ with respect to $\pi_{P_1}(Y)$. Note also that $\pi_{P_1}(Y)$ is a set of 6 points in LGP in $\PP^3$ since $Z$ is in LGP. Thus $\pi_{P_1}(Q)$ avoids the Weddle surface of $\pi_{P_1}(Y)$, so $\pi_{\pi_{P_1}(Q)} (\pi_{P_1}(Y))$ is a set of 6 points not on a conic in $\PP^2$. 

But 
\[
\pi_{\pi_{P_1}(Q)} (\pi_{P_1}(Y)) = \pi_\lambda(Y) = \pi_{\pi_Q(P_1)}(\pi_Q(Y))
\]
and by assumption these six points in $\PP^2$ do lie on a conic. This is a contradiction.
\end{proof}

\begin{proof}[Proof of Proposition \ref{thequad}]
By Lemma \ref{nocon} in a general projection $Z'$ of $Z$ no point is a vertex for some quadric $\mathcal Q$ through $Z'$. By Lemma \ref{l. 2b points on a quadric} there exists a curve $C$ of degree $b$ through $Z'$, which splits in a union of lines, such that the intersection $\mathcal Q\cap C$ is finite. The claim follows by B\'ezout's Theorem.    
\end{proof}

It follows that in order to determine the structure of $(b,2)$-geprofi sets in LGP, with $b\geq 4$ we can limit ourselves to study sets of points $Z$ in $\PP^4$ whose general projection to $\PP^3$ lies in a quadric surface. Clearly we want that $Z$ is nondegenerate, so in particular it cannot lie in a surface of degree 2, by the LGP property.

The problem of sets in $\PP^4$ whose general projection lies in a quadric is interesting by itself, and we will treat it as a self-standing problem.

\begin{example}\label{b2gen} It is a trivial consequence of the previous discussion that any set of $2b\leq 8$ points in linear general position in $\PP^4$ are geprofi. 

Namely for any such set a general projection $Z'$ to a hyperplane is contained in an irreducible quadric (this is clear for $2b=8$, while for $2b<8$ it follows since we can find quadrics through $Z'$ which split in totally different pairs of planes,  so we take a general linear combination). Then we can apply Lemma \ref{l. 2b points on a quadric}.
\end{example}

The following result, proved in  \cite[Lemma 2.1]{maroscia}, and its consequences on the Hilbert function of  points in LGP, will be often useful in the rest of the paper.

\begin{theorem}[Maroscia]\label{Maroscia} Let $Z$ be a  finite set of points in $\PP^n$, with $h$-vector $h=(1, h_1, \dots, h_s)$, where $s$ is the smallest integer such that $Z$ is separated in degree $s$, so that $h_s>0$. Assume that  $h_i=n-t$ for some $0< i<s$ and for some $t>0$. Then $Z$ contains a subset of at least $n-t+2$ points contained in a hyperplane of $\PP^n$.
In particular $Z$ is not in LGP.
    \end{theorem}

  Maroscia's Theorem can be considered as a preliminary result to the classical Castelnuovo's Lemma, that we will need below, and we state here for the reader's convenience.

\begin{lemma}  [Castelnuovo's Lemma] \label{Cast} Let $Z$ be a set of $q$ points in LGP in $\PP^n$. Assume $q\geq 2n+3$ and assume that the Hilbert function of $Z$ satisfies $H_Z(2)\leq 2n+1$. Then  in fact $H_Z(2) = 2n+1$ and $Z$ lies in a rational normal curve in $\PP^n$.
\end{lemma}

\begin{lemma} \label{pencil of quadrics}
Let $Z$ be a set of $q\geq 9$ points in $\mathbb P^4$ in LGP. Then the general projection of $Z$ does not lie on a pencil of quadrics.
\end{lemma}
\begin{proof} It is enough to prove the statement for $q=9$.  Suppose, for contradiction, that a general projection of $Z$ does lie on a pencil of quadrics. Let $Z_1$ and $Z_2$ be two disjoint sets of three points of $Z$. Let $H_1$ and $H_2$ be the  planes  spanned by $Z_1$ and $Z_2$ and let $P$ be their point of intersection. We project from $P$. So $\pi_P(Z_1)$ and $\pi_P(Z_2)$ are two sets of three collinear points in $\mathbb P^3$. By semicontinuity, the set $\pi_P(Z)$ is  also contained in (at least) a pencil of quadrics.  Thus, in particular, it is contained in a quadric cone $\mathcal Q$. By the LGP assumption the line spanned by $\pi_P(Z_1)$ and the line spanned by $\pi_P(Z_2)$ are skew, hence $\mathcal Q$ splits into the union of two planes. Therefore, $\pi_P(Z)$ has at least five points in a plane and $Z$ has five points in a 3-dimensional linear space. This contradicts the LGP assumption.       
\end{proof}

\begin{remark} \label{classical weddle}
    The classical Weddle surface is the  surface of degree 4 arising as the 2-Weddle locus (see \cite[Definition 2.1]{POLITUS1}) of a set $Z$ of six points in $\PP^3$ in linear general position (LGP). Equivalently, this is the set of points in $\PP^3$ from which the projection to $\PP^2$ is a set of points lying on a conic. Section 2.3 of \cite{POLITUS1} explores different behavior that can arise for Weddle surfaces of six points when the LGP condition is relaxed. In the same paper, Example 2.8 studies the 2-Weddle locus for 10 general points in $\PP^4$, and in Remark 8.34 it is shown that for a general set $Z$ of 8 points in $\PP^4$, the locus of points $P$ from which the  projection $\pi_P(Z)$ to $\PP^3$ is a complete intersection of three quadrics is a union of lines plus a ``somewhat mysterious arithmetically Cohen-Macaulay curve of degree 7 in $\PP^4$".

    Now we ask a related question. For a set of 9 points in $\PP^4$ in LGP, what are the possible 2-Weddle loci? 
    In Lemma \ref{pencil of quadrics} we showed that the general projection of $Z$ does not lie on a pencil of quadrics, so the 2-Weddle locus represents the points of projection from which we get the image lying in a pencil of quadrics in $\PP^3$. Of course any line joining two of the points must lie on the 2-Weddle locus, but we can do more. 

    Let $Z = \{ P_1,\dots,P_9\}$ and let $L_i$ be the linear form dual to $P_i$ for each $i$. As in \cite[Chapter~2]{POLITUS1} we use Macaulay duality. By Maroscia's theorem \cite{maroscia} and the LGP assumption, the $h$-vector of $Z$ must be $(1,4,4)$ so in particular $Z$ imposes independent conditions on quadrics, so we have
    \[
\dim [I_Z]_2 = 15-9 = 6 = \dim [R/(L_1^2,\dots,L_9^2)]_2.
    \]
    For a general point $P$ and its dual general linear form $L$, we want to compute  
    \[
    \dim [I_Z \cap I_P^2]_2 = \dim [R/(L_1^2,\dots,L_9^2,L)]_2.
    \]
    Consider the exact sequence
    \[
 [R/(L_1^2,\dots,L_9^2)]_1 \stackrel{\times L}{\longrightarrow} [R/(L_1^2,\dots,L_9^2)]_2 \rightarrow 
    [R/(L_1^2,\dots,L_9^2,L)]_2 \rightarrow 0.
    \]

The first vector space has dimension 5 while the second has dimension 6, so multiplication by $L$ is represented by a $6 \times 5$ matrix of linear forms, and the ideal of maximal minors defines a proper subvariety of $\PP^4$ by Lemma \ref{pencil of quadrics}. Thus the codimension of the 2-Weddle locus is either 2 or 1.

The expected codimension of the 2-Weddle locus is 2. This expected codimension is achieved for a general set of 9 points:  indeed, by direct computation, the codimension for a random set of 9 points is the expected one, so by semicontinuity it happens also for a general set of 9 points. 
Note that then the degree of the 2-Weddle locus is necessarily $\binom{6}{2} = 15$ (\cite[Lemma 1.4]{geom-inv}).

It remains to show that it can also happen that the 2-Weddle locus is a hypersurface. Let $Z$ be a set of 9 points on a rational normal curve $C$ in $\PP^4$. Let $W$ be the secant variety of $C$, which is a hypersurface of degree 3, see \cite[Theorem 4.3]{Dale84}. The projection $\pi$ from a general point of $W$ sends $C$ to a nodal curve $\pi(C) \subset \PP^3$ of degree 4. Thus its arithmetic genus is 1 and it is a complete intersection of quadrics. Hence $\pi(Z)$ lies on a pencil of quadrics. Thus the 2-Weddle locus of $Z$ is a hypersurface.
\end{remark}

We will focus then on the case where the general projection of $Z$ lies in a unique quadric surface.

\begin{example}\label{on2}
    Let $C\subset\PP^4$ be a rational normal (quartic) curve. Then by Riemann-Roch a general projection $C'$ of $C$ in $\PP^3$ lies in a quadric surface $\mathcal Q$. Since $C'$ is not linearly normal,  it is not a complete intersection, and thus $\mathcal Q$ cannot be a cone by \cite[Exercise V.2.9]{hartshorne}. Hence $C'$ is a curve of type either $(1,3)$ or $(3,1)$ in a smooth quadric surface. 

    It is clear that any set of  points in $C$ projects generically in $\PP^3$ to points contained in a smooth quadric surface. Thus, for any even $q=2b>4$ a general set $Z$ of $q$ points in $C$ determines an example of a $(b,2)$-geprofi set. It is well known that $Z$, as any subsets of $C$, must be in LGP.

  Finally, suppose $Z$ is a set of $4k+2$ points on $C$, $k \geq 2$. Then $Z$ has $h$-vector $(1,4,4,\dots,4,1)$. No subset of $\geq 5$ points on $C$ is degenerate, so by Maroscia's Theorem \ref{Maroscia} $Z$ has the Cayley-Bacharach Property, and hence by the main result of \cite{DGO} $Z$ is arithmetically Gorenstein.
    
\end{example}

\subsection{The (5,2) case} We focus first on the case $b=5$. Since 10 general points in $\PP^3$ are contained in no quadrics, it is clear that a general set of 10 points in $\PP^4$ is not a $(5,2)$-geprofi set. Nevertheless we will see that there are sets in LGP,  not lying on an elliptic quintic, whose general projection is contained in one (necessarily irreducible) quadric. 
\smallskip

\begin{remark} \label{h2=4} Let $Z$ be a set of 10 points in LGP in $\PP^4$. Then by Maroscia's Theorem \ref{Maroscia} the $h$-vector of $Z$ starts with $1,4,h_2,\dots$ where $h_2$ cannot be smaller than 4. 

Assume that $h_2=4$, so that the $h$-vector of $Z$ is $(1,4,4,1)$, and $Z$ is contained in 6 independent quadric hypersurfaces. Fix $P\in Z$ and consider the set $Z_0=Z\setminus \{P\}$. The $h$-vector of $Z_0$ is $(1,4,4)$ since otherwise, by Maroscia's Theorem \ref{Maroscia}, $Z_0$ would not be in LGP. In particular, any quadric containing $Z_0$ also contains $Z$. It is elementary that for a general $Q\in\PP^4$ there exists a quadric cone with vertex at $Q$ containing $Z_0$. Thus there exists a quadric cone with vertex $Q$ that contains $Z$, i.e., a general projection of $Z$ lies in a quadric surface.
 It follows by Proposition~\ref{thequad} that $Z$ is a $(5,2)$-geprofi set.
\end{remark}

\begin{example}\label{no4}
    Consider a set of $6$ general points in a hyperplane $H\subset\PP^4$.  Take $4$ general quadrics in $\PP^4$ vanishing at these $6$ points. They intersect in $16$ points. Take the set $Z$ of $10$ points residual to the original $6$. The $h$-vector of $Z$ is $(1,4,4,1)$ by construction, and the set is in linear general position (computed by Singular). Thus $Z$ is a $(5,2)$-geprofi, by Remark~\ref{h2=4}.

    Notice that $Z$ is  not contained in  rational normal quartic curves, since such a quartic would be in the base locus of the system of quadrics vanishing in the points. But the four chosen quadrics in the construction of $Z$ intersect in $16$ points only (computed by Singular). 
\end{example}

\begin{remark}
    Example \ref{no4} shows that Castelnuovo's Lemma \ref{Cast} does not  hold for $2n+2$ points in $\PP^n$.
\end{remark}

We have found two examples of $(5,2)$-geprofi sets in LGP: sets of 10 points in a rational  normal (quartic) curve,  and  sets constructed as in Example \ref{no4} by liaison. 
We observe that both of these examples have $h$-vector $(1,4,4,1)$, and both are arithmetically Gorenstein sets (see Example \ref{on2} -- the argument for the situation in Example \ref{no4} is identical.)

\begin{remark} 
Let $C$ be a nondegenerate smooth curve of degree 5 and genus 1 in $\PP^4$. Let $Z \in |\mathcal O_C(2)|$. Since $C$ is ACM, there is a quadric hypersurface that cuts out $Z$, so $\dim [I_Z]_2 = \dim [I_C]_2 + 1 = 6$. Then using the methods of \S\ref{LGPcurves} below, one can check that $Z$ is a nontrivial $(5,2)$-geprofi set. Nevertheless, it is not hard to see that it participates in a link using four quadric hypersurfaces, with  residual a set of 6 points with generic Hilbert function inside a  hyperplane, and no three points on a line. We do not know if every set $Z$ constructed as in Example \ref{no4} lies on such a curve $C$.
\end{remark}

To  classify all $(5,2)$-geprofi in LGP with $h$-vector $(1,4,4,\dots)$ we need the following result.

\begin{lemma}\label{B4ic}
    Let $Z$ be a set of $q\geq 9$ points in LGP in $\PP^4$, whose $h$-vector starts with $(1,4,4,\dots)$. Call $B$ the base locus of the system of quadrics through $Z$. Then either $B$ is finite, or $B$ is an irreducible rational quartic  curve. 
\end{lemma}
\begin{proof}
   We start by observing that no reducible quadric can contain $Z$, by the LGP assumption. In particular, $B$ cannot have dimension $2$. It also follows that for any hyperplane $L$, the system of quadrics through $Z$ restricts to a system in $L$ of (affine) dimension $6$ of quadric surfaces in $L$, containing $L\cap Z$. 

   Assume that $B$ contains an irreducible curve $Y$. If $Y$ spans $\PP^4$, take a  general  set $W$ of more than  $4 \cdot \deg(Y)$ points in $Y$. $W$ is in LGP and  lies in a 6-dimensional system of quadrics, so, by Castelnuovo's Lemma \ref{Cast}, $W$ is contained in an irreducible rational quartic~$\Gamma$. But no quartic curve can meet $Y$ in more than $4 \cdot \deg(Y)$ points (e.g., by B\'ezout on a general plane projection of $Y$ and $\Gamma$). Thus $\Gamma=Y$. Since $Y$ is cut out by quadrics and the system of quadrics through $Y$ is 6-dimensional, then $B=Y$. In  particular, $Y$ contains $Z$.

   Assume $Y$ is degenerate. No irreducible curves that span $\PP^3$ lie in a 6-dimensional system of quadrics  in $\mathbb P^3$. Thus $Y$ is a plane curve.

   Assume $\deg(Y)\geq 2$, so that $Y$ spans a plane $\pi$. Fix a hyperplane $L$ containing $\pi$, fix a general point $A_1\in\pi$, and fix general points $A_2,\dots,A_5\in L$.  Since $Z$ is in LGP, no reducible quadric in $\PP^4$ can contain $Z$. The linear system of quadrics through $Z$ has affine dimension~$6$, and no members of the system contain $L$. Thus the system cuts on $L$ a linear system of quadric surfaces of affine dimension~$6$. It follows that 
 there is a quadric surface $S$ in $L$, cut by a quadric hypersurface passing through $Z$, which contains $\{A_1,\dots,A_5\}$, and necessarily $Y\subset S$.  Since $S$ contains $Y$ and $A_1$, then it contains $\pi$ (even when $\deg(Y)=2$).  Hence $S$ must split in a union of two planes. This is impossible because no plane in $L$ contains four general points $A_2,\dots,A_5$.

   So assume that $Y$ is a line. Necessarily no more than two points of $Z$ lie in $Y$. Take a hyperplane $L$  containing $Y$ and two more points $P_1,P_2\in Z$. Fix general lines $\ell_1,\ell_2$ in $L$, through $P_1,P_2$ resp. Fix general points $A_2,A_3\in\ell_1$ and $A_4,A_5\in\ell_2$. As above, the linear system of quadrics through $Z$  cuts on $L$ a linear system of quadric surfaces of affine dimension $6$. Thus
   there is a pencil of quadric surfaces  in $L$, cut by quadric hypersurfaces passing through $Z$, all of which contain $\{A_2,\dots,A_5\}$, and necessarily also $Y, P_1,P_2$.  The quadrics of the pencil contain the three pairwise skew lines $Y,\ell_1,\ell_2$. This cannot happen in $\PP^3$, so we get a contradiction.
\end{proof}

Now we consider the case where the $h$-vector of $Z$ is $(1,4,5)$. We prove below that this case cannot occur.

\begin{proposition}\label{145} There is no $(5,2)$-geprofi set in LGP in $\PP^4$, whose $h$-vector is $(1,4,5)$.
\end{proposition}
\begin{proof}
  Suppose that  $Z$ were a set of 10 points in LGP in $\PP^4$, whose $h$-vector is $(1,4,5)$, and whose general projection to $\PP^3$ lies in a quadric surface. We will show that this leads to a contradiction. Then for a general point $P\in\PP^4$  there is a quadric cone $\mathcal{Q}_P$ containing  $Z$, with vertex in $P$. When $P$ moves generically, we find a family $\mathcal F$ of quadric cones, whose vertices are dense in $\PP^4$. 

    Let $Z'\subset Z$ be a set of 9 points (necessarily in LGP, in particular nondegenerate) and call $P_1$ the point of $Z$ not in $Z'$.  Since the $h$-vector of $Z'$ is necessarily $(1,4,4)$, we see that  the linear  system $L_{Z'}$ of quadrics through $Z'$, viewed as a linear space in the $\PP^{14}$ that parametrizes quadrics in $\PP^4$, has dimension 5, while the corresponding system $L_Z$ for $Z$ has dimension 4. By Lemma \ref{pencil of quadrics}, for a general $P\in\PP^4$ there exists a unique quadric  through $Z'$ which is a cone with vertex $P$. Since $L_{Z'}$ is 5-dimensional, the general quadric through $Z'$  is smooth. Then the space $L_{Z'}$ meets the locus of quadric cones (which is defined by the determinant of a $5\times 5$ matrix of variables) in a hypersurface $W$ of degree 5 in $L_{Z'}$. 
    
    The family $\mathcal F$ is parametrized by a subvariety $F$ of $W\cap L_Z$ inside $L_{Z'}$. 
    Assume first that $F$ has dimension 4. We have $F \subset W \cap L_Z \subset L_Z$ so all of these varieties have dimension 4. 
    Then $F$ is dense in $L_Z$, so by semicontinuity the elements of $L_Z$ are cones whose vertices cover $\PP^4$. But, by Bertini's Theorem, the singularities of a general element of a linear system are contained in the base locus of the system, which cannot be $\PP^4$.

    Thus either $F$ is not dense in $L_Z$, or there exists a proper (nonlinear) subfamily $\mathcal F'$ of $\mathcal F$ of quadric cones whose vertices are dense in $\PP^4$.

    After replacing $\mathcal F$ with a smaller family, we may always assume that $\dim(F)\leq 3$. Since the vertices of the elements of $\mathcal F$ are dense in $\PP^4$, and $Z$ cannot lie in a reducible quadric by the LGP assumption, then necessarily $\dim F=3$ and the general element of $\mathcal F$ is a quadric cone of rank 3.
    Thus, for a general $P\in \PP^4$ there exists a quadric cone of rank $\leq 3$ containing $Z$, with vertex in a line passing through $P$. In other words the projection from $P$ sends $Z$ to a set of points in LGP in $\PP^3$ which are contained in a quadric cone. 

    Now we obtain the final contradiction. Consider 6 points $P_2,\dots,P_7$ of $Z$ and let $P$ be the intersection of the two planes $\pi,\pi'$ spanned by $P_2,P_3,P_4$ and $P_5,P_6,P_7$ resp. The  LGP property guarantees that no line joining two points of $Z$ passes through $P$ (else we have five points in a hyperplane), and $\pi\cup\pi'$ spans $\PP^4$. By semicontinuity, the projection of $Z$ from $P$ lies in a quadric cone $Y$ in $\PP^3$. But there are two skew lines in $\PP^3$, projection of $\pi,\pi'$, each containing three points of $Z$. Thus $Y$ contains two skew lines. This means $Y$ splits in a union of two planes. This contradicts the LGP property of $Z$.    
\end{proof}

\begin{corollary}\label{c. LGP + 1441}
    Let $Z$ be a set of 10 points in LGP in $\PP^4$. Then $Z$ is (5,2)-geprofi if and only if $h_Z=(1,4,4,1).$ 
\end{corollary}
\begin{proof}
    From Proposition \ref{145} a (5,2)-geprofi set in LGP has $h$-vector $(1,4,4,1). $ Conversely, if $Z$ has $h$-vector $(1,4,4,1)$ then $\dim[{I_Z}]_2=6$ and a general double point imposes at most 5 conditions.
\end{proof}

\subsection{The $(b\geq 6,2)$ case}

The classification of  $(b,2)$-geprofi sets for $b\geq 6$  is easier than the classification of the $(5,2)$-case, because we have more than 10 points and we can invoke Castelnuovo's Lemma \ref{Cast}. We will also use the following fact.

\begin{theorem} [{\cite[Theorem 2.5]{GMR}}] \label{GMR fact}
Let $s$ be a positive integer and let $X$ be a projective variety containing more than $s$ points. Assume the Hilbert function of $X$ is $H_X(i)$, $i \geq 0$. Then there is a finite subset $V \subset X$ such that the Hilbert function of $V$ is
\[
H_V(i) = \min \{ H_X(i), s \}.
\]

\end{theorem}

\begin{corollary}\label{b>5}
    Any $(b,2)$-geprofi set in LGP in $\PP^4$, with $b\geq 6$, is contained in a rational normal quartic.
\end{corollary}
\begin{proof}
    Let the $h$-vector of $Z$ start with $1,4,h_2, \dots$.  Since $|Z| \geq 12$, it follows from Maroscia's theorem and LGP that $h_2 \geq 4$. If $h_2=4$, then $H_Z(2)=9$. Again since $Z$ has at least 12 points, it lies in a rational normal quartic by  Castelnuovo's Lemma.

    Now assume that $h_2 \geq 5$. By Theorem \ref{GMR fact} it follows  that  $Z$ contains a  subset $Z'$ of 10 points with $h$-vector $1,4,5$.  This produces a contradiction because the general projection of $Z$ lies in a quadric, which means that the general projection of $Z'$ lies in a quadric and $Z'$ is in LGP, a situation excluded by Proposition \ref{145}.
\end{proof}

Thus we can summarize the classification of all $(b,2)$-geprofi sets in LGP.

\begin{theorem} \label{2b LGP class}
    Let $Z$ be a set of $2b$ points in LGP in $\PP^4$. 
  Then $Z$ is a $(b,2)$-geprofi set if and only if one of the following occurs.
    \begin{itemize}
\item[(i)]$b=3$ (in fact such $Z$ is always a trivial geprofi set).

\item[(ii)] $b=4$. 

\item[(iii)] $b=5$, $Z$ is arithmetically Gorenstein, and either it lies on a rational quartic curve, or it is linked by a system of quadrics to a degenerate scheme of length $6$.

\item[(iv)]$b\geq 6$ and $Z$ lies on a rational normal curve.  (If $b$ is odd then $Z$ is arithmetically Gorenstein.)
    \end{itemize}
\end{theorem}

\begin{proof}
    The case $b\leq 4$ is contained in Example \ref{b2gen}.

    When $b=5$ we know by Proposition \ref{145} and Remark \ref{h2=4} that the $h$-vector of $Z$ starts with $1,4,4,\dots$. Call $B$ the base locus of the system of quadrics through $Z$. By Lemma~\ref{B4ic}, either $B$ is finite or $B$ is a rational normal quartic. When $B$ is finite then a general system of 4 quadrics through $Z$ links $Z$ to a set $W$ of 6 points. The liaison properties imply that the $h$-vector of $W$ starts with $1,3,\dots$, so that $W$ is degenerate. On the other hand, by Example~\ref{no4}, if $W$ is a general set of $6$ points in a hyperplane of $\PP^4$, a general liaison of $W$ with a system of quadrics produces a $(5,2)$-geprofi set. If $B$ is a rational normal quartic, then every set of $10$ points on $B$ is in LGP, and it is $(5,2)$-geprofi by Remark \ref{h2=4}. 

    When $b>5$ then $Z$ lies in a rational normal quartic by Corollary \ref{b>5}. Such examples exist by Remark \ref{h2=4}.  Since $Z$ lies on a rational normal curve, its $h$-vector is either $(1,4,\dots,4,3)$ or $(1,4,\dots,4,1)$, and it has the Cayley-Bacharach property (in fact the stronger Uniform Position Property). So in the latter case it is arithmetically Gorenstein by \cite{DGO}.
\end{proof}

\begin{remark}
We notice that the same argument used in Remark \ref{h2=4} proves that, regardless of the LGP property, a general projection of any set $Z$ of $2b$ points whose $h$-vector starts with $1,4,4,\dots$ is contained in a quadric surface. Indeed, suppose $Z$ has $h$-vector starting $(1,4,4,\dots)$ but does not necessarily have LGP. 
    By Theorem \ref{GMR fact}, we {\it can find} a subset $Z' \subset Z$  so that the set $Y = Z \backslash Z'$ of 9 points has $h$-vector $(1,4,4)$. Notice that $[I_Y]_2 = [I_Z]_2$. Then the projection of $Y$ to a hyperplane lies on a quadric, and as before then there is a quadric cone with vertex at a general point that contains $Z$.

\end{remark}


\section{Geprofi sets on curves in $\PP^4$} \label{LGPcurves}

In this section we give a useful construction of geprofi sets in $\PP^4$. This will have two consequences. First, it will be a primary component of the proof of a numerical classification of geprofi sets in $\PP^4$ not necessarily in LGP (Theorem \ref{numerical classification}), which is one of the main results of this paper. Second, it fills in roughly half of the possibilities working toward a numerical classification of geprofi sets in $\PP^4$ in LGP (see Remark \ref{LGP table}). Again, we point out the significant contrast with $\PP^3$, where conjecturally only a set of four general points gives a  geproci set in LGP, and where there are no $(a,b)$-geproci sets on an irreducible curve in $\PP^3$ of degree $\max\{a,b\}$ or less, as was shown in \cite{POLITUS3}.

We start with the following useful observation.

\begin{proposition}\label{3x3-7}
Let $Z=\{P_1,\dots,P_7\}$ be a nondegenerate set of $7$ points in $\PP^4$. Then a general projection of $Z$ to $\PP^3$ lies in no twisted cubics.
\end{proposition}

\begin{proof} Notice that if $Z$ has 3 collinear points or 4 coplanar points, then $Z$ cannot project to a twisted cubic, because all finite subsets of twisted cubics are in LGP.
Thus we may assume that a general projection of $Z$ to a hyperplane  is in LGP in that hyperplane, since otherwise we are immediately done. In fact, even projection from $P_1$ to the hyperplane sends $P_2,\dots,P_7$ to a set of six points in LGP in the hyperplane.

Let $L$ be a general hyperplane in $\PP^4$ and let $L'$ be a general hyperplane containing $P_1$. Let $Q'$ be a general point of $L$; by construction $Q'$ is a general point of $\PP^4$. Let $M = L \cap L'$. 
Call $\pi_1$ the projection from $P_1$ to $L$ and $\pi'$ the projection from $Q'$ to $L'$.

Suppose, for the sake of contradiction, that a general projection of $Z$ to $L'$ lies on a twisted cubic in $L'$. So we can assume that $\pi'(P_1),\dots,\pi'(P_7)$ lie on a twisted cubic in $L'$. Notice that $\pi'(P_1) = P_1$. 

Notice also that $\pi_1 \circ \pi'$ sends points of $\PP^4$, away from the line $\lambda_{P_1Q'}$ joining $P_1$ and $Q'$, to $M$. Also, $\pi_1 \circ \pi' = \pi' \circ \pi_1$ away from $\lambda_{P_1Q'}$.

Now consider the points $P_2,\dots,P_7$, which are all off of $\lambda_{P_1Q'}$. 
The points $\pi_1(P_2),\dots,\pi_1(P_7)$ are in linear general position in $L$. Since $Q'$ is general, it is not in the Weddle quartic of 
$\pi_1(P_2),\dots,\pi_1(P_7)$ (see Remark \ref{gen proj 6 pts}), so the projections $\pi'(\pi_1(P_2)),\dots,\pi'(\pi_1(P_7))$ lie in no conics of $M$. 

On the other hand, for $P_2,\dots,P_7$ we have 
$\pi_1 \circ \pi' = \pi' \circ \pi_1$. But we have assumed that $\pi'(P_1), \pi'(P_2),\dots,\pi'(P_7)$ lie on a twisted cubic in $L'$ and that $\pi'(P_1) = P_1$. But now projecting from a point of a twisted cubic sends the rest of the points to points on a conic. Thus
\[
\{ \pi'(\pi_1(P_2)),\dots,\pi'(\pi_1(P_7)) \} = \{ \pi_1 (\pi' (P_2), \dots, \pi_1 (\pi'(P_7)) \}
\]
do lie on a conic, giving us our contradiction. 
\end{proof}

\begin{corollary}\label{no33LGP}
    There are no $(3,3)$-geprofi sets in LGP in $\PP^4$.
\end{corollary}
\begin{proof}
    A $(3,3)$-geprofi set in $\PP^4$  in LGP  would project generically to a twisted (irreducible) cubic curve in $\PP^3$. 
\end{proof}

Motivated by Proposition \ref{3x3-7} we propose the following problem.

\begin{problem}
    Let $Z \subset \PP^4$ be a set of points such that the general projection  $\pi(Z)$ to $\PP^3$ lies on a smooth curve $C'$ of degree $b\ge 4$. Find an (effective if possible) bound $N(b)$ so that if $|Z| \geq N(b)$ then $Z$ itself lies on a smooth curve $C$ of degree $b$ in $\PP^4$. If we also require $\pi(C) = C'$ then of course $C'$ must fail to be linearly normal. Is this latter requirement automatic when $|Z| \geq N(b)$? Or are there situations where $\pi(Z)$ lies on a linearly normal curve of degree $b$ but $Z$ itself does not lie on a curve of degree $b$?
\end{problem}  

Our next result gives a construction of nontrivial $(b,d)$-geprofi sets on curves of $\PP^4$.

\begin{proposition} \label{g=b-4}
Let $Z' \subset C' \subset \PP^3$, where $Z'$ is a set of points and $C'$ is an irreducible nondegenerate curve of degree~$b\ge 4$. Assume that $Z'$ is cut out on $C'$ by a surface of degree $d\ge 2$. Assume:

\begin{itemize}

\item[(i)] $C'$ is smooth;

\item[(ii)] $C'$ is not linearly normal;

\item[(iii)] $h^1(\mathcal I_{C'}(d)) = 0$.

\end{itemize}

\noindent Then there exists $Z \subset  \PP^4$ lying on a smooth curve $C$ of degree $b$, such that $Z$ is a nontrivial $(b,d)$-geprofi set.
\end{proposition}

\begin{proof}

Note that we start with $b=4$, because  all irreducible nondegenerate curves of degree 3 are linearly normal (they are twisted cubics).

Let $V = H^0(\mathcal O_{C'}(1))$. By assumption, $\dim V \geq 5$. Thus the linear system $|{\mathcal O}_{C'}(1)|$ embeds $C'$ as a curve $C \subset \PP^4$ (possibly after projecting down to $\PP^4$). The projection sending $C$ to $C'$ amounts to studying a 4-dimensional subspace of $V$. Since $C' \cong C$, we have the corresponding set $Z$ on $C$. By semicontinuity, for a general projection $\pi : \PP^4 \rightarrow \PP^3$, we have $h^1(\mathcal I_{\pi(C)}(d)) = 0$. Thus $\pi (Z) $ is cut out in $\pi(C)$ by a surface of degree $d$, and hence $Z$ is geprofi. The last part of the proposition is immediate from Proposition \ref{triv irred}. 
\end{proof}

\begin{remark}\label{r. GPL}
    By the main theorem of \cite{GLP}, condition (iii) of Proposition \ref{g=b-4} holds for all smooth curves $C' \subset \PP^3$ when $d \geq b -2$. Furthermore, as long as $C'$ is not a rational curve with a $(b -1)$-secant line, then it also holds for $d = b-3$.
\end{remark}

As a preparation for the next result we recall the following fact.

\begin{theorem}[\cite{BE-inv}] \label{BE max rk}
If $b \geq \frac{3g+12}{4}$, $g \geq 0$, there exists a smooth connected curve in $\PP^3$ of degree $b$ and genus $g$ such that
for all $t \geq 0$, the restriction map
\[
H^0(\mathcal O_{\PP^n}(t)) \rightarrow H^0(\mathcal O_C(t))
\]
has maximal rank.    
\end{theorem}

Note that the numerical condition in this theorem is equivalent to $g \leq \frac{4}{3} b -4$. The next result shows that if we fix $b$, then there is a ``reasonable" lower bound $M(b)$ for $d$ so that we can find a smooth curve $C$ in $\PP^4$ of degree $b$ and a divisor $Z \in |\mathcal O_C (d)|$ for any $d \geq M(b)$ so that $Z$ is geprofi.

\begin{corollary} \label{find geprofi by curves}

Let $C'$ be a general irreducible smooth curve of degree $b$ and genus $g=b-4$ in $\mathbb P^3$.  Fix $d>2$.  Assume 

\begin{itemize}

\item[(i)] $d > - \frac{10}{b}+2$ 

\item[(ii)] $b(d-1) +5 \leq \binom{d+3}{3}$.

\end{itemize}

Then $C'$ is the projection of a curve $C$ in $\PP^4$, and a general degree $d$ hypersurface section $Z$ of $C$ is $(b,d)$-geprofi. By generality, $Z$ is a nontrivial geprofi set in linear general position.
\end{corollary}

\begin{proof}

By Theorem \ref{BE max rk}, there exists a smooth curve of degree $d$ and genus $g$ in $\PP^3$ having maximal rank. 

The assumption that $g=b-4$ guarantees that $C'$ will not be linearly normal, by Riemann-Roch.

Condition (i) is equivalent to $2g-2 < bd$. This implies, in particular, that $h^1(\mathcal O_C(d)) = 0$ since $bd$ is greater than the degree of the canonical divisor. Then the Riemann-Roch formula gives
\[
h^0(\mathcal O_C(d)) = bd - g + 1 + h^1(\mathcal O_C(d)) = bd - (b-4) + 1 = b(d-1) +5.
\]
Now consider the exact sequence of sheaves
\[
0 \rightarrow \mathcal I_{C'}(d) \rightarrow \mathcal O_{\PP^3}(d) \rightarrow \mathcal O_{C'}(d) \rightarrow 0.
\]
In the long exact sequence in cohomology, condition (ii) and Theorem \ref{BE max rk} give that 
\[
H^0(\mathcal O_{\PP^3}(d)) \rightarrow H^0(\mathcal O_{C'}(d))
\]
 is surjective, so $h^1(\mathcal I_{C'}(d)) = 0$. We then apply Proposition \ref{g=b-4} to get the result.
\end{proof}

\begin{example}
Now we will use different properties of a smooth rational quartic curve to show that $(4,d)$-geprofi point sets exist for any $d \geq 3$. Indeed, if $C'$ is a smooth rational quartic curve in $\PP^3$ then the restriction map
\[
H^0(\mathcal O_{\PP^3}(d)) \rightarrow H^0(\mathcal O_{C'}(d))
\]
is surjective for all $d \geq 2$ (see Remark \ref{r. GPL}). But $C'$ is the projection from $\PP^4$ of a rational normal curve $C$. Thus any element $Z$ of the linear system $|\mathcal O_C(d)|$ is in LGP and has the property that for a general projection $\pi$, $\pi(Z)$ is cut out by a surface of degree $d$, so $Z$ is $(4,d)$-geprofi. This is just the first case of Corollary \ref{find geprofi by curves}.

By the same token, we use Corollary \ref{find geprofi by curves} to get $(b,d)$-geprofi sets for 
\[
\begin{array}{lcccc}
b=5,\ d \geq 2, \\
b=6,\ d \geq 3, \\
b=7,\ d \geq 3, \\
b=8,\ d \geq 4, \\
b=9,\ d \geq 4, \\
b=10,\ d \geq 4, \\
b=11, \ d \geq 5.
\end{array}
\]
\end{example}

\begin{remark}\label{r. b or d 2}
It is clear that $(2,d)$-geprofi sets cannot be in linear general position for $d \geq 3$ since a reduced curve of degree 2 can only be a plane conic (possibly reducible) or a pair of skew lines, and either way we violate linear general position.  On the other hand, Theorem \ref{2b LGP class} shows that $(b,2)$-geprofi sets in LGP exist for all $b \geq 2$, but for $b=2$ they are degenerate.

Summarizing, if $Z$ is a $(3,d)$-geprofi set with $d \geq 3$ (which must be trivial), then

    \begin{itemize}
        \item[a)] $Z$ is not in LGP.

        \item[b)] One of the following two possibilities holds.

        \begin{itemize}
        \item[i)] $Z$ is equally distributed in three pairwise disjoint lines.
        
        \item[ii)] $Z$ has $d$ points on a line $\ell$ and $2d$ points of a (possibly reducible) conic $C$, such that the plane of $C$ is disjoint from $\ell$; moreover, if $C$ is reducible then $Z$ has exactly  $d$ points on each line of $C$.
        \end{itemize}
    \end{itemize}
\end{remark}

\begin{remark}\label{r. (3,d) geprofi}
Theorem \ref{2b LGP class} shows that $(3,2)$-geprofi sets exist (although such a set in $\PP^4$ is a trivial geprofi set). 
We now claim that a $(3,d)$-geprofi set in linear general position does not exist for any $d \geq 3$. Indeed, suppose such a set $Z$ did exist. For $Z$ to be in linear general position, for a general projection $\pi$ we must have that $\pi(Z)$ is the intersection of a twisted cubic curve and a surface of degree $d$. In particular, $\pi(Z)$ is a set of $\geq 9$ points on a twisted cubic curve.   Since such a curve lies on (more than) a pencil of quadrics, this contradicts Lemma \ref{pencil of quadrics}.
\end{remark}

We can now complete the numerical classification of geprofi sets in $\PP^4$ (not necessarily in LGP), which is analogous to the numerical classification of geproci sets in $\PP^3$ given in \cite[Theorem 4.6]{POLITUS1}.

\begin{theorem} \label{numerical classification}
A nondegenerate, nontrivial $(b,d)$-geprofi set exists in $\PP^4$ if and only if  $b \geq 4$ and $d \geq 2$.  (In particular, there is no nontrivial $(3,d)$-geprofi set in $\PP^4$ for any $d$.)
\end{theorem}

\begin{proof}
We will use the fact that a general projection of a set of points $Z$ in $\PP^4$ is collinear if and only if $Z$ itself is collinear, and similarly if we replace ``collinear" with ``coplanar".

For the case $b=2$ see Remark \ref{r. b or d 2}.

Now set $b=3$ and $d=2$. Let $Z$ be a set of 6  points in $\PP^4$ in LGP. As in Theorem \ref{2b LGP class}, $Z$ is $(3,2)$-geprofi, but it is trivial. 

Similarly, let $b=3$ and $d \geq 3$. Suppose that a nontrivial $(3,d)$-geprofi set $Z$ exists in $\PP^4$. The general projection $\pi(Z)$ of $Z$ lies on a curve $C'$ of degree $3$, and the image is cut by a surface of degree $d$. We rule out that $C'$ is a twisted cubic by Proposition \ref{3x3-7}. So $C'$ is either the nondegenerate union of a line and a conic or else the nondegenerate union of three skew lines. Respectively, $Z$ must consist of $d$ points on any such line, and $2d$ points on any such conic. Since the projection is general and $d \geq 3$, $Z$ must lie similarly on lines and conics. But then taking three points at a time we see that $Z$ is a trivial geprofi set.

Now fix any $b \geq 4$ and $d \geq 2$. If $b \geq d+2$ then we have by Theorem \ref{t.extend}  that a nontrivial $(b,d)$-geprofi  set exists (not necessarily in LGP). If 
\[
b \leq \frac{1}{d-1} \left [ \binom{d+3}{3} -5 \right ]
\]
then we have by Corollary \ref{find geprofi by curves} that a nontrivial $(b,d)$-geprofi set exists (even in LGP). But it is elementary to check that 
\[
\frac{1}{d-1} \left [ \binom{d+3}{3} -5 \right ] > d+2
\]
for $d \geq 2$ so these two results cover all remaining cases.    
\end{proof}

\begin{remark} \label{LGP table}
We have seen that the LGP condition restricts the possible existence of geprofi sets in $\PP^4$. Nevertheless, it is interesting that so many examples in LGP do exist, since the only known geproci set in LGP in $\mathbb P^3$ consists of four points  (cf. \cite[Question 0.2]{POLITUS1}).

Let us put together the information we have on nondegenerate, nontrivial $(b,d)$-geprofi sets in $\PP^4$ in linear general position. In Table \ref{tab LGP}, $E$ means there is a nondegenerate $(b,d)$-geprofi set in LGP in $\PP^4$ for the given $b$ and $d$, $X$ means it does not exist, and  ``?'' means we do not yet know. An arrow $\rightarrow$ or $\downarrow$ means that the pattern ($E$ or $X$ or ``?") continues forever after that point. Of course the cases $b=1$ and $d=1$ are trivial (by linear general position) and are left out. Note also that in this table we are interested more in the LGP property than in the nontriviality property; in particular, the case $(3,2)$ exists in LGP but it is a trivial geprofi set. Also, the case $(2,2)$ is degenerate and so is left out for that reason.

\begin{table}[h!]
    \centering
\[
\begin{array}{r|ccccccccccccccc}
 &d:  \\
\text{LGP} & 2&3&4&5&6&7&8&9&10& \dots \\ \hline
b: \ \ 2 & X & X & X & X & X & X & X & X & X &  \rightarrow \\
3 & E & X & X & X & X & X & X & X & X  & \rightarrow \\
4 & E & E & E & E & E & E & E & E & E & \rightarrow   \\
5 & E & E & E & E & E & E & E & E & E & \rightarrow   \\
6 &  E & E & E & E & E & E & E & E & E & \rightarrow  \\
7 & E & E & E & E & E & E & E & E & E & \rightarrow   \\
8 & E & ?& E & E & E & E & E & E & E  & \rightarrow   \\
9 & E & ? & E & E & E & E & E & E & E  & \rightarrow  \\
10 & E &? & E & E & E & E & E & E &  E & \rightarrow   \\
11 & E &? &?& E & E & E & E & E &  E & \rightarrow   \\
12 & E & ? & ? &  E & E & E & E & E &  E & \rightarrow   \\
13 & E & ? & ? &  ? & E & E & E & E &  E & \rightarrow   \\
 & \downarrow & \downarrow & \downarrow & \downarrow  &  &  & &  & 
\end{array}
\]
    \caption{The map of geprofi sets in $\PP^4$ in LGP.}
    \label{tab LGP}
\end{table}
\end{remark}

From 
Theorem \ref{t.extend} 
we know that 
if we fix $d$ then we can produce a nontrivial $(b,d)$-geprofi set for any $b \geq d+2$. However, the set we produce has many collinearities, and it is of interest to ask the same question for points in LGP. 

\begin{question} \label{final question} If we fix $d$, can we find a  value $M(d)$ so that for $b > M(d)$  there is no  $(b,d)$-geprofi set  in LGP? In particular, is it true that for fixed $d$ and for all $b \gg d$ there is no $(b,d)$-geprofi set   in LGP in $\PP^4$?
\end{question}

Here we have produced many examples of $(b,d)$-geprofi sets in LGP, but only for $b$ relatively small compared to $d$.
This suggests that the answer to the latter question is ``yes" (in contrast to the situation of Theorem~\ref{numerical classification} where we do not assume LGP). We point out that a similar situation occurs for geproci sets in $\PP^3$: there are none we know of with five or more points in LGP (cf. \cite{POLITUS1}).


\vspace{.3in}

\paragraph*{\bf Acknowledgement.}
This work was initiated while all authors participated in the Research Group ``Intersections in Projective Spaces'' hosted by the BIRS program in Kelowna in summer 2023. We thank BIRS for its generous support and in particular Chad Davis for making our stay in British Columbia so pleasant.

We also thank the ``Welcome to Poland'' program of NAWA,
grant no. PPI/WTP/2022/1/ 00063 which supported the Conference ``Lefschetz Properties in Algebra, Geometry, Topology and Combinatorics'' in Krak\'ow in summer 2024, during which our work has been completed.

\bibliographystyle{abbrv}

\end{document}